\documentclass[10pt,reqno]{amsart}

\usepackage[utf8]{inputenc}
\usepackage[english]{babel}
\usepackage[dvipsnames]{xcolor}
\usepackage{amssymb}
\usepackage{hyperref}
\usepackage{mathrsfs}
\usepackage{mathtools}
\mathtoolsset{showonlyrefs}

\usepackage{xparse}% http://ctan.org/pkg/xparse
\NewDocumentCommand{\ceil}{s O{} m}{%
  \IfBooleanTF{#1} % starred
    {\left\lceil#3\right\rceil} % \ceil*[..]{..}
    {#2\lceil#3#2\rceil} % \ceil[..]{..}
}

\usepackage{color}
\usepackage{enumerate}

\newtheorem{theorem}{Theorem}[section]
\newtheorem{lemma}[theorem]{Lemma}
\newtheorem{corollary}[theorem]{Corollary}
\newtheorem{proposition}[theorem]{Proposition}

\theoremstyle{remark}
\newtheorem{remark}[theorem]{\bf{Remark}}

\theoremstyle{definition}
\newtheorem{assumption}[theorem]{Assumption}
\newtheorem{example}[theorem]{Example}
\newtheorem{definition}[theorem]{Definition}

\newcommand\cbrk{\text{$]$\kern-.15em$]$}}
\newcommand\opar{\text{\,\raise.2ex\hbox{${\scriptstyle
|}$}\kern-.34em$($}}
\newcommand\cpar{\text{$)$\kern-.34em\raise.2ex\hbox{${\scriptstyle |}$}}\,}

\newcommand{\aint}{-\hspace{-0.38cm}\int}
\newcommand{\aaint}{-\hspace{-0.31cm}\int}

\newcommand\bC{\mathbb{C}}

\newcommand\bE{\mathbb{E}}

\newcommand\bH{\mathbb{H}}

\newcommand\bL{\mathbb{L}}

\newcommand\bN{\mathbb{N}}

\newcommand\bP{\mathbb{P}}
\newcommand\bQ{\mathbb{Q}}
\newcommand\bR{\mathbb{R}}

\newcommand\bZ{\mathbb{Z}}

\newcommand\cF{\mathcal{F}}
\newcommand\cG{\mathcal{G}}

\newcommand\cM{\mathcal{M}}

\newcommand\cO{\mathcal{O}}

\newcommand\cS{\mathcal{S}}

\newcommand\nn{\nonumber}

\newcommand{\mysection}[1]{\section{#1}
\setcounter{equation}{0}}

\newcommand{\Ccinf}{C_{c}^{\infty}}
\newcommand{\R}{\mathbb{R}}

\begin{document}

\title[PDE with nonlocal operators having slowly varying symbols]
{An $L_q(L_p)$-theory for space-time non-local equations generated by L\'evy processes with low intensity of small jumps}

\author{Jaehoon Kang}
\address{Stochastic Analysis and Application Research Center, Korea Advanced Institute of Science and Technology, 291 Daehak-ro, Yuseong-gu, Daejeon, 34141, Republic of Korea}
\email{jaehoon.kang@kaist.ac.kr}

\author{Daehan Park}
\address{Stochastic Analysis and Application Research Center, Korea Advanced Institute of Science and Technology, 291 Daehak-ro, Yuseong-gu, Daejeon, 34141, Republic of Korea} 
\email{daehanpark@kaist.ac.kr}

\thanks{The authors were  supported by the National Research Foundation of Korea(NRF) grant funded by the Korea government(MSIT) (No. NRF-2019R1A5A1028324).}



\begin{abstract}
We investigate an $L_{q}(L_{p})$-regularity ($1<p,q<\infty$) theory for space-time nonlocal equations of the type $\partial^{\alpha}_{t}u = \mathcal{L}u +f$.  
Here, $\partial^{\alpha}_{t}$ is the Caputo fractional derivative of order $\alpha\in(0,1)$ and $\mathcal{L}$ is an integro-differential operator
$$
\mathcal{L}u(x) = \int_{\mathbb{R}^{d}} \left( u(x)-u(x+y) -\nabla u (x) \cdot y \mathbf{1}_{|y|\leq 1} \right) j_{d}(|y|)dy
$$
which is the infinitesimal generator of an isotropic unimodal L\'evy process. We assume that the jump kernel $j_{d}(r)$ is comparable to $r^{-d} \ell(r^{-1})$, where $\ell$ is a continuous function satisfying 
$$
C_{1}\left(\frac{R}{r}\right)^{\delta_{1}} \leq \frac{\ell(R)}{\ell(r)} \leq C_{2} \left( \frac{R}{r} \right)^{\delta_{2}} \quad \text{for}\;\; \,1\leq r\leq R<\infty,
$$
where $0\leq \delta_{1}\leq \delta_{2}<2$. Hence, $\ell$ can be slowly varying at infinity. Our result covers $\mathcal{L}$ whose Fourier multiplier $\Psi(\xi)$ satisfies
$\Psi(\xi)\asymp -\log{(1+|\xi|^{\beta})}$ for $\beta \in (0,2]$ and $\Psi(\xi) \asymp-(\log(1+|\xi|^{\beta/4}))^{2}$ for $\beta\in(0,2)$ by taking $\ell(r) \asymp 1$ and $\ell(r) \asymp \log{(1+r^{\beta})}$ for $r\geq1$ respectively. In this article, we use the Calder\'on-Zygmund approach and function space theory for operators having slowly varying symbols.
\end{abstract}

\maketitle

\tableofcontents

\mysection{Introduction}
\label{Main result sec}

Equations with space or time nonlocal operators are used to model natural phenomena in various area of science (see e.g. \cite{metzler1999anomalous,bouchaud1990,fogedby1994}). For example, time-fractional heat equation $\partial^{\alpha}_{t}=\Delta u $ ($\alpha\in(0,1)$) describes subdiffusive aspect of anomalous diffusion caused by particle sticking and trapping effect. Also, when we give relativistic correction to Laplacian, then it becomes relativistic Hamiltonian $-(\sqrt{-\Delta +m^{2}}-m)$ which is a nonlocal operator generated by relativistic stable process.

In this article, we consider equations with space-time nonlocal operator
\begin{equation}\label{eqn 09.03.16:02}
\partial^{\alpha}_{t}u=\mathcal{L}u +f \quad t>0,
\end{equation}
where $\partial^{\alpha}_{t}$ is the Caputo fractional derivative of order $\alpha\in(0,1)$, and $\mathcal{L}$ is an integro-differential operator whose jump kernel can have low intensity of small jumps. Precisely speaking, we have the following representation for $\mathcal{L}$:
\begin{align}\label{op_rep}
\mathcal{L}u(x) &= \int_{\bR^{d}} \left( u(x)-u(x+y) -\nabla u(x) \cdot y \mathbf{1}_{|y|\leq 1} \right) j_{d}(|y|) dy,
\end{align}
where $j_{d}(r)$ is decreasing in $r$ and comparable to $r^{-d}\ell(r^{-1})$, and $\ell$ is a positive continuous function which satisfies
$$
C_{1}\left(\frac{R}{r}\right)^{\delta_{1}} \leq \frac{\ell(R)}{\ell(r)} \leq C_{2} \left( \frac{R}{r} \right)^{\delta_{2}} \quad \text{for}\;\; \,1\leq r\leq R<\infty,
$$
where $0\leq \delta_{1}\leq \delta_{2}<2$. We remark that $\ell$ can be slowly varying at infinity (see Definition \ref{d:sv} and Assumption \ref{ass bernstein} for detail) due to the possible choice of $\delta_{1}=0$. Therefore, $j_{d}$ can have low intensity of small jumps. Examples of function $\ell$ covered by our result are
\begin{gather}\label{eqn 09.25.13:01}
\ell(r) = r^{\beta}, \quad \ell(r)=r^{\beta}/(r^{\beta}+1), \quad \ell(r)=\log{(1+r^{\beta})} \quad \text{for} \quad \beta \in (0,2). \nonumber
\end{gather}

It is known that $\mathcal{L}$ can be considered as a linear operator with symbol $-\psi$, where $\psi$ is the L\'evy-Khintchine exponent of L\'evy process whose jump kernel is $j_d(|y|)$. Our spatial nonlocal operator is the infinitesimal generator of pure jump L\'evy process whose transition density function $p_{d}(t,x)$ is the fundamental solution to parabolic equation $u_{t} = \mathcal{L}u$. Under our assumptions, the  heat kernel for jump processes which generates $\mathcal{L}$ has different type of estimation for small time and large time (see Section \ref{sec 11.09.16:14}). Moreover, since we are dealing with equations with time-fractional derivatives, the fundamental solution to \eqref{eqn 09.03.16:02} is the transition density of time changed process by an inverse subordinator (see Section 4). This implies that our estimation for the fundamental solution needs more exquisute analysis.

In the literature, equations with non-local operators in  time variables have been widely studied. See for example, \cite{clement1992global,chen2020time,dong2021p,kim17timefractionalpde,P13,zacher2005maximal}. Regarding (parabolic) equations with spatial nonlocal operator $u_{t}=Lu +f$, where $L$ is of the form
$$
Lu(t,x) = \int_{\bR^{d}} \left( u(t,x+y)-u(t,x)-\nabla_{x}u(t,x) \cdot y \chi^{(\sigma)}(y) \right) J(t,x,y) dy,
$$
an $L_{p}$-estimation of solution was introduced in \cite{mik1992}. Here, $\chi^{(\sigma)}$ is a function depending on $\sigma \in (0,2)$ and $J(t,x,y) = a(t,x,y)|y|^{-d-\sigma}$, where $a(x,y)$ is homogeneous of order zero and sufficiently smooth in $y$. Most of the studies focus on $J$ which generalizes $a(t,x,y)$. See e.g. \cite{kim2019L,mikulevivcius2017p,mikulevivcius2019cauchy,zhang2013p}. 
Quite recently, \cite{dong2021nonlocal} proved an $L_{p}$-estimation of solution to equations
$$
\partial^{\alpha}_{t}u = Lu +f,
$$
where $\alpha\in(0,1]$ (i.e. the result also covers parabolic case) and $J$ is comparable to $|y|^{-d-\sigma}$ uniformly in $(t,x)$ and H\"older continuous in $x$ uniformly in $(t,y)$.

An $L_{q}(L_{p})$-regularity result was introduced in \cite[Theorem 8.7]{P13} and \cite{kim2020nonlocal}. The result in \cite[Theorem 8.7]{P13} deals with abstract parabolic Volterra equations of the form
$$
u(t) + \int_{0}^{t} a(t-s) Au(s) ds = f(t),
$$
where $a$ is locally integrable function and $A$ is densely defined closed operator on $L_{p}$. The class of $A$ is general and it covers operators $-\phi(-\Delta)$ for Bernstein functions $\phi$. In \cite{kim2020nonlocal}, %the $L_{q}(L_{p})$-regularity for the solutions to 
the following equation is studied:
\begin{align}\label{eqn 09.25.19:08}
\partial^{\alpha}_t u=-\phi(-\Delta)u+f,\quad t >0;\quad u(0,\cdot)=u_0,
\end{align}
where $\alpha\in(0,1)$ and $\phi$ is a Bernstein function satisfying
\begin{align}\label{phi_lsc}
c\left(\frac{R}{r}\right)^{\delta_0}\le\frac{\phi(R)}{\phi(r)},\quad\text{for all}\;\;0<r\le R<\infty \quad (c>0, \quad \delta_{0}\in(0,1]).
\end{align}
The authors in \cite{kim2020nonlocal} only used elementary analysis, based on estimation of the heat kernel $p_{d}(t,x)$
\begin{align}\label{hke_lsc}
|p_d(t,x)| \le C\left(\phi^{-1}(t^{-1})^{d/2}\wedge \frac{t\phi(|x|^{-2})}{|x|^d}\right),
\end{align}
(see \cite{kim2013parabolic,kim2020nonlocal}) which highly depends on scaling condition \eqref{phi_lsc}.

Although \cite[Theorem 8.7]{P13} and \cite{kim2020nonlocal} can cover a large class of operators, it is hard to see if we can apply the results to an integro-differential operator $\mathcal{L}$ of the form \eqref{op_rep}. To apply 
\cite[Theorem 8.7]{P13} or \cite{kim2020nonlocal}, $\mathcal{L}$ should satisfy the condition in \cite[Theorem 8.7]{P13} or the symbol $\psi$ of $\mathcal{L}$ is written as $\psi(\xi)=\phi(|\xi|^2)$ for a Bernstein function $\phi$.
However, it seems that it is complicated to check whether the operator $\mathcal{L}$ (or its symbol $\psi$) satisfies the above conditions. In fact, \cite[Theorem 8.7]{P13} requires a comprehensive background in abstract harmonic analysis to check the conditions therein even for $A=-\phi(-\Delta)$ (see e.g. \cite[Section 3.2]{kim2022spde}). Moreover, it is difficult to check that $\psi$  satisfies $\psi(\xi)=\phi(|\xi|^2)$ for a Bernstein function $\phi$. To the best of the author's knowledge, the only known direct relation between $\psi$ and $j_{d}$ is
\begin{align*}
\psi(\xi)=\int_{\R^d}(1-\cos(\xi \cdot x))j_d(|x|)dx,
\end{align*}
which follows from L\'evy-Khintchine formula, and it is hard to obtain closed form of $\psi$ from the above integral in general. Even if we know closed form of $\psi$, it is another problem to check that $\psi$ can be written as $\psi(\xi)=\phi(|\xi|^2)$. Thus, it is not easy to see if we can apply \cite[Theorem 8.7]{P13} or \cite{kim2020nonlocal} to the $\mathcal{L}$ in \eqref{op_rep}. % when we do not know its symbol. 
In addition, $\phi$ should satisfy scaling condition \eqref{phi_lsc} to apply \cite{kim2020nonlocal}.

Motivated by the above observation, this paper aims to find general conditions on $\mathcal{L}$ that give $L_{q}(L_{p})$-regularity of solutions to \eqref{eqn 09.03.16:02}. This extends the results in \cite[Theorem 8.7]{P13} and \cite{kim2020nonlocal}  in the following two aspects: 
\begin{itemize}
\item[(1)] Our results require elementary analysis and cover more general operators. i.e., the symbol $\psi$ of operator needs not satisfy $\psi(\xi)=\phi(|\xi|^2)$ for a Bernstein function $\phi$. Indeed, we do not need neither closed form of $\psi$ nor smoothness of $\psi$.

\item[(2)] Also, our result does not need scaling condition \eqref{phi_lsc} even for $\mathcal{L}=-\phi(-\Delta)$.
\end{itemize}

\smallskip

\noindent To achieve this goal, we give assumptions on $j_d$ instead of symbol of operator $\mathcal{L}$ in \eqref{op_rep}. Indeed, our assumption on $j_d$ is quite general so that our results cover all operators of the form $-\phi(-\Delta)$ for $\phi(r)=\log(1+r^{\beta})$ with $\beta\in(0,1]$,  and for a Bernstein function $\phi$ satisfying \eqref{phi_lsc} and ${\phi(R)}/{\phi(r)} \leq c (R/r)^{\delta_{0}'}$ for all  $0<r\leq R<\infty$,
where $0<\delta_0\le\delta_0'<1$ (see \cite{BGR14a,KM12,KSV14}).

It is worth mentioning that our approach and \cite{kim2020nonlocal,P13} have their own advantages. If the symbol $\psi$ of an operator $\mathcal{L}$ is explicitly given as $\psi(\xi)=\phi(|\xi|^2)$ for a Bernstein function $\phi$, then  
\cite{kim2020nonlocal,P13} is more accessible since one needs not to know asymptotic behavior of jump kernel of the operator $\mathcal{L}$. On the contrary, if an operator $\mathcal{L}$ is given as \eqref{op_rep}, then our approach is more accessible since we do not have to obtain exact value of symbol of $\mathcal{L}$. 

In this paper, we investigate maximal $L_{q}(L_{p})$-regularity of solutions to \eqref{eqn 09.03.16:02}, where $\mathcal{L}$ is an integro-differential operator generated by an isotropic unimodal pure jump L\'evy process $X$ in $\R^d$. As mentioned above, we assume that the jump kernel of the L\'evy process $X$ is comparable to $|x|^{-d}\ell(|x|^{-1})$ and consider weak scaling conditions on $\ell$ which cover the case that $\ell$ is a slowly varying function at infinity. In this case, we may have $\sup_x p_d(t,x)=\infty$ and thus we cannot have \eqref{hke_lsc} for our heat kernel.

To obtain our main results, we follow the standard approach in harmonic analysis. Precisely, we control the sharp function of derivative of solution in terms of $L_{\infty}$-norm of free term $f$, and then use the Fefferman-Stein theorem and the Calder\'on-Zygmund theorem. Main difficulty arises here since $X$ is an isotropic unimodal L\'evy process which is a more general process than the one in \cite{kim2020nonlocal}. Here we give a short description. Unlike \cite{kim2020nonlocal}, we cannot expect global scaling condition like \eqref{phi_lsc} to underlying functions related to $X$. Hence, $p_{d}(t,x)$ has different form of estimation comparing to \eqref{hke_lsc} and thus our proof is much more involved. Indeed, if $\limsup_{r\to\infty}\ell(r)=\infty$, then $\ell$ gives the borderline between near and off diagonal estimates. Since the scaling function for parabolic cube is $\psi$, and the two functions $\ell$ and $\psi$ may not be comparable, we need more delicate argument.

We finish the introduction with some notations. We use $``:="$ or $``=:"$ to denote a definition. The symbol $\bN$ denotes the set of positive integers and $\bN_0:=\bN\cup\{0\}$. Also, we use $\bZ$ to denote the set of integers. As usual $\bR^d$ stands for the Euclidean space of points $x=(x^1,\dots,x^d)$. We set
$$
B_r(x):=\{y\in \bR: |x-y|<r\}, \quad \bR_+^{d+1} := \{(t,x)\in\bR^{d+1} : t>0 \}.
$$
For $i=1,\ldots,d$,
multi-indices $\sigma=(\sigma_{1},\ldots,\sigma_{d})$,
and functions $u(t,x)$ we set
$$
\partial_{x^{i}}u=\frac{\partial u}{\partial x^{i}}=D_{i}u,\quad D^{\sigma}u=D_{1}^{\sigma_1}\cdots D_{d}^{\sigma_d}u,\quad|\sigma|=\sigma_{1}+\cdots+\sigma_{d}.
$$
We also use the notation $D_{x}^{m}$ for arbitrary partial derivatives of
order $m$ with respect to $x$.
For an open set $\cO$ in $\bR^{d}$ or $\bR^{d+1}$, $C_c^\infty(\cO)$ denotes the set of infinitely differentiable functions with compact support in $\cO$. By
$\cS=\cS(\bR^d)$ we denote the class of Schwartz functions on $\bR^d$.
For $p> 1$, by $L_{p}$ we denote the set
of complex-valued Lebesgue measurable functions $u$ on $\R^{d}$ satisfying
\[
\left\Vert u\right\Vert _{L_{p}}:=\left(\int_{\R^{d}}|u(x)|^{p}dx\right)^{1/p}<\infty.
\]
Generally, for a given measure space $(M,\mathcal{M},\mu)$, $L_{p}(M,\cM,\mu;F)$
denotes the space of all $F$-valued $\mathcal{M}^{\mu}$-measurable functions
$u$ so that
\[
\left\Vert u\right\Vert _{L_{p}(M,\cM,\mu;F)}:=\left(\int_{M}\left\Vert u(x)\right\Vert _{F}^{p}\mu(dx)\right)^{1/p}<\infty,
\]
where $\mathcal{M}^{\mu}$ denotes the completion of $\cM$ with respect to the measure $\mu$.
If there is no confusion for the given measure and $\sigma$-algebra, we usually omit the measure and the $\sigma$-algebra.
We denote $a\wedge b := \min\{a,b\}$ and $a\vee b:=\max\{a,b\}$. By $\cF_{d}$ and $\cF^{-1}_{d}$ we denote the $d$-dimensional Fourier transform and the inverse Fourier transform respectively, i.e.
$$
\cF_{d}(f)(\xi):=\hat{f}(\xi):=\int_{\bR^d} e^{-i\xi\cdot x} f(x)dx, \quad \cF^{-1}_{d}(f)(\xi):=\frac{1}{(2\pi)^d}\int_{\bR^d} e^{i\xi\cdot x} f(x)dx.
$$
For two real-valued functions $f,g$ defined on a set $A$, we write $f\asymp g$ on $A$ if there is a constant $c>1$ such that $c^{-1}f(b)\leq g(b)\leq cf(b)$ for all $b\in A$. Finally if we write $C=C(\dots)$, this means that the constant $C$ depends only on what are in the parentheses. The constant $C$ can differ from line to line.

\medskip

\section{Main results}
\label{Main result sec}
In this section, we introduce our main results. We first present our spatial nonlocal operator $\mathcal{L}$. Let $X=(X_t, t\ge0)$ be a  L\'evy process on $\bR^d$ with L\'evy-Khintchine exponent $\psi$. Then,
\begin{align*}
\bE e^{i \xi \cdot X_t} =\int_{\R^d}e^{i \xi \cdot x}p_d(t,dx)=e^{-t\psi(\xi)},
\end{align*}
where $p_d(t, dx)$ is the transition probability of $X_{t}$. If $X$ is a pure jump symmetric L\'evy process with L\'evy measure $j_{d}$, then $\psi$ is of the form
\begin{align*}
\psi(\xi)=\psi_X(\xi)=\int_{\R^d}(1-\cos(\xi \cdot x))j_d(dx),
\end{align*}
where $\int_{\R^d} (1\wedge |x|^2)j_d(dx)<\infty$.

A measure $\mu(dx)$ is isotropic unimodal if it is absolutely continuous on $\R^d\setminus\{0\}$ with a radial and radially decreasing density. A L\'evy process $X$ is isotropic unimodal if $ p_d(t, dx)$ is isotropic unimodal for all $t > 0$. This is equivalent to the condition that the L\'evy measure $j_d(dx)$ of $X$ is isotropic unimodal if $X$ is pure jump L\'evy process (see \cite{W83}).

Throughout this paper, we always assume that $X$ is a pure jump isotropic unimodal L\'evy process with the L\'evy-Khintchine exponent $\psi$. 
With a slight abuse of notation, we will use the notations $\psi(|x|) = \psi(x)$  and $j_d(dx)=j_d(x)dx=j_d(|x|)dx$ for $x \in \R^d$.

For $f\in \cS(\bR^{d})$, define a linear operator $\mathcal{L}$ as
$$
\mathcal{L}f(x) = \int_{\bR^{d}} \left( f(x)-f(x+y)-\nabla f(x) \cdot y \mathbf{1}_{|y| \leq 1} \right) j_{d}(y)dy.
$$
Due to the L\'evy-Khinchine representation, $\psi$ is continuous and negative definite (see \cite[Theorem 1.1.5]{farkaspsi}). Thus, by \cite[Proposition 2.1.1]{farkaspsi} we can understand $\mathcal{L}$ as the infinitesimal generator of $X$, and nonlocal operator with Fourier multiplier $-\psi(|\xi|)$. Precisely speaking, for $f\in \cS(\bR^{d})$,  we have the following relation
\begin{equation}\label{eqn 06.28.11:33}
\mathcal{L}f(x) = \lim_{t \downarrow 0} \frac{ \bE f(x+X_t) -f(x)}{t} = \cF^{-1}(-\psi(|\xi|)\cF(f)(\xi))(x).
\end{equation}
In this context, we also use notations $\mathcal{L}_{\psi}$ or $\mathcal{L}_{X}$ instead of $\mathcal{L}$ in this article.

One of well-known examples of isotropic unimodal L\'evy process is subordinate Brownian motion $Y=(Y_{t},t\geq0)$ which is defined by $Y_t:=B_{S_{t}}$. Here $B=(B_{t},t\geq0)$ is  a $d$-dimensional Brownian motion and $S=(S_{t}, t\geq0)$ is a subordinator (i.e., $1$-dimensional increasing L\'evy process) independent of $B$. It is known that there is a Bernstein function $\phi: \bR_+ \to \bR_+$ (i.e. $(-1)^{n}\phi^{(n)}\le0$ for all $n\in\bN$, where $\phi^{(n)}$ is the $n$-th derivative of $\phi$) satisfying $\bE[e^{-\lambda S_t}]=e^{-t\phi(\lambda)}$. If $\phi(0+)=0$, then $\phi$ has the following representation
\begin{equation*}
\phi(\lambda)=b\lambda + \int_{(0,\infty)} (1-e^{-\lambda t})\mu(dt),
\end{equation*}
where $b\geq 0$ and $\mu$ is a measure  satisfying  $\int_{(0,\infty)} (1\wedge t) \mu(dt)<\infty$. $\mu$ is called the L\'evy measure of $\phi$ (see \cite{SSV12}).

It is well-known that the L\'evy-Khintchine exponent of $Y$ is $\xi\mapsto\phi(|\xi|^2)$ and the L\'evy measure of $Y$ has the density $J_d(x)=J_d(|x|)$, where
\begin{equation}\label{eqn 04.28.16:34}
J_{d}(r)=\int_{(0,\infty)} (4\pi t)^{-d/2} e^{-r^2/(4t)} \mu(dt).
\end{equation}
Hence, $\mathcal{L}_{Y}$ also has representation \eqref{eqn 06.28.11:33} with $J_{d}$ and $\phi(|\cdot|^{2})$ in place of $j_{d}$ and $\psi$. For example, by taking $\phi(\lambda)=\lambda^{\beta/2}$ ($\beta\in (0,2)$), we obtain the fractional Laplacian $\Delta^{\beta/2}=-(-\Delta)^{\beta/2}$, which is the infinitesimal generator of a rotationally symmetric $\beta$-stable process in $\bR^d$. 

In order to describe the regularity of solution, we introduce Sobolev space related to the operator $\mathcal{L}$.
For $\gamma\in\bR$, and $u\in \cS(\bR^{d})$, define linear operators 
\begin{equation*}
(-\mathcal{L})^{\gamma/2}=(-\mathcal{L}_{\psi})^{\gamma/2}, \quad (1-\mathcal{L})^{\gamma/2}=(1-\mathcal{L}_{\psi})^{\gamma/2}
\end{equation*}
as follows
\begin{align*}
\cF\{(-\mathcal{L})^{\gamma/2}u\}= (\psi(|\xi|))^{\gamma/2}\cF(u)(\xi),\quad
\cF\{(1-\mathcal{L})^{\gamma/2}u\}= (1+\psi(|\xi|))^{\gamma/2}\cF(u)(\xi).
\end{align*}
For $1<p<\infty$, let $H_p^{\psi,\gamma}$ be the closure of $\cS(\bR^{d})$ under the norm
\begin{equation*}
\|u\|_{H_p^{\psi,\gamma}}:=\|\cF^{-1}\{\left(1+\psi(|\cdot|)\right)^{\gamma/2}\cF(u)(\cdot)\}\|_{L_p}<\infty.
\end{equation*}
Then from the definition of $H^{\psi,\gamma}_{p}$ the operator $(1-\mathcal{L})^{\gamma/2}$ can be extended from $\cS(\bR^{d})$ to $L_{p}$. Throughout this article, we use the same notation $(1-\mathcal{L})^{\gamma/2}$ for this extension. For more information, see e.g. \cite{farkaspsi}.  Also note that if $\psi(|\xi|)=|\xi|^{2}$, then $H^{\psi,\gamma}_{p}$ is a standard Bessel potential space $H^{\gamma}_{p}$ and $H^{\psi,0}_{p}=L_{p}$ due to the definition.

The following lemma is a collection of useful properties of $H^{\psi,\gamma}_{p}$.
\begin{lemma}\label{H_p^phi,gamma space}
Let $1<p<\infty$ and let $\gamma\in\bR$.

(i) The space $H_p^{\psi,\gamma}$ is a Banach space.

(ii) For any $\mu\in\bR$, the map $(1-\mathcal{L})^{\mu/2}$ is an isometry from $H^{\psi,\gamma}_{p}$ to $H^{\psi,\gamma-\mu}_{p}$.

(iii) If $\mu>0$, then we have continuous embeddings $H_p^{\psi,\gamma+\mu}\subset H_p^{\psi,\gamma}$ in the sense that
\begin{equation*}
\|u\|_{H_p^{\psi,\gamma}}\leq C \|u\|_{H_p^{\psi,\gamma+\mu}},
\end{equation*}
where the constant $C$ is independent of $u$. 

(iv) For any $u\in H^{\psi,\gamma+2}_{p}$, we have
\begin{equation}\label{eqn 03.25.15:03}
\left(\|u\|_{H^{\psi,\gamma}_p}+\|\mathcal{L}u\|_{H^{\psi,\gamma}_p}\right) \asymp  \|u\|_{H_p^{\psi,\gamma+2}}.
\end{equation}

\end{lemma}

\begin{proof}
The first and second assertions are direct consequences of the definition. Recall that $\psi$ is a continuous negative definite function. Hence, the third assertion comes from \cite[Theorem 2.3.1]{farkaspsi}. Finally, the last assertion can be obtained by using the second assertion and  \cite[Theorem 2.2.7]{farkaspsi}.
\end{proof}

Now we introduce our non-local operator in time variable and related definitions.
For $\alpha>0$ and $\varphi\in L_{1}((0,T))$, the Riemann-Liouville fractional integral
of the order $\alpha$ is defined as
$$
I_{t}^{\alpha}\varphi:=\frac{1}{\Gamma(\alpha)}\int_{0}^{t}(t-s)^{\alpha-1}\varphi(s)ds, \quad 0\leq t\leq T.
$$
We also define  $I^0\varphi:=\varphi$.  Take $n\in \bN$ such that  $\alpha \in [n-1, n)$. 
If  $\varphi(t)$ is  $(n-1)$-times  differentiable and $\left(\frac{d}{dt}\right)^{n-1} I_t^{n-\alpha}  \varphi$ is absolutely continuous on $[0,T]$, then
the Riemann-Liouville fractional derivative $D_{t}^{\alpha}$ and the Caputo fractional derivative $\partial_{t}^{\alpha}$ are defined as
\begin{equation}
                          \label{eqn 4.15}
D_{t}^{\alpha}\varphi:=\left(\frac{d}{dt}\right)^{n}\left(I_{t}^{n-\alpha}\varphi\right),
\end{equation}
and
\begin{align*}
\partial_{t}^{\alpha}\varphi= D_{t}^{\alpha} \left(\varphi(t)-\sum_{k=0}^{n-1}\frac{t^{k}}{k!}\varphi^{(k)}(0)\right).
\end{align*}
 Using  Fubini's theorem, we see that for any $\alpha,\beta\geq 0$,
 \begin{equation}
                                                          \label{eqn 4.15.3}
I^{\alpha}_tI^{\beta}_t \varphi=I^{\alpha+\beta}_t \varphi, \quad
\text{$(a.e.)$} \,\, t\leq T.
\end{equation}
Note that $D^{\alpha}_t\varphi=\partial^{\alpha}_t \varphi$ if $\varphi(0)=\varphi'(0)=\cdots=\varphi^{(n-1)}(0)=0$.
By \eqref{eqn 4.15.3} and \eqref{eqn 4.15},  if  $\alpha,\beta\geq 0$,
\begin{equation*}
                \label{eqn 4.20.1}
D^{\alpha}_tD^{\beta}_t=D^{\alpha+\beta}_t, \quad D^{\alpha}_t I_{t}^{\beta} \varphi=
D_{t}^{\alpha-\beta}\varphi,
\end{equation*}
where we define $D_t^{a}\varphi:=I_t^{-a}\varphi$ for $a<0$.
Also if 
$\varphi(0)=\varphi^{(1)}(0)=\cdots = \varphi^{(n-1)}(0)=0$ 
then by definition of $\partial^{\alpha}_{t}$,
\begin{equation*}
I^{\alpha}_{t}\partial^{\alpha}_{t}u=I^{\alpha}_{t}D^{\alpha}_{t}u=u.
\end{equation*}

For $p,q\in(1,\infty), \gamma\in\bR$ and $T<\infty$, we denote
\begin{equation*}
\bH_{q,p}^{\psi,\gamma}(T):=L_q\left((0,T); H_p^{\psi,\gamma}\right), \qquad \bL_{q,p}(T):=\bH_{q,p}^{\psi,0}(T).
\end{equation*}
We write $u\in C_{p}^{\alpha,\infty}([0,T]\times\bR^d)$ if $D^m_x u, \partial^{\alpha}_t D^m_x u \in C([0,T];L_p)$ for any $m\in \bN_{0}$. Also $C_{p}^{\infty}(\bR^d)=C^{\infty}_{p}$ denotes the set of functions $u_0=u_0(x)$ such that $D^m_xu_0\in L_p$ for any $m\in \bN_{0}$. 

\begin{definition}
 \label{defn defining}
Let $\alpha\in(0,1)$, $1<p,q<\infty$, $\gamma\in\bR$, and $T<\infty$.

(i) We write $u\in {\bH_{q,p}^{\alpha,\psi,\gamma+2}}(T)$ if there exists a sequence $u_n\in C^{\alpha,\infty}_{p}([0,T]\times \bR^d)$ satisfying
\begin{equation*}
\|u-u_n\|_{\bH_{q,p}^{\psi,\gamma+2}(T)} \to 0 \quad \text{and} \quad \|\partial^{\alpha}_{t}u_{n} - \partial^{\alpha}_{t}u_{m} \|_{\mathbb{H}^{\psi,\gamma}_{q,p}(T)} \to 0
\end{equation*}
as $n,m \to \infty$. We call this sequence $u_n$ a defining sequence of $u$, and we define
\begin{equation*}
\partial^{\alpha}_t  u= \lim_{n \to \infty} \partial^{\alpha}_t u_n \text{ in } \bH_{q,p}^{\psi,\gamma}(T).
\end{equation*}
The norm in $\bH_{q,p}^{\alpha,\psi,\gamma+2}(T)$ is naturally given by
\begin{equation*}
\|u\|_{{\bH_{q,p}^{\alpha,\psi,\gamma+2}}(T)}=\|u\|_{\bH_{q,p}^{\psi,\gamma+2}(T)}+\|\partial^{\alpha}_{t}u\|_{\bH_{q,p}^{\psi,\gamma}(T)}.
\end{equation*}

(ii) We write $u\in{\bH_{q,p,0}^{\alpha,\psi,\gamma+2}(T)}$, if there is a defining sequence $u_{n}$ of $u$ such that $u_{n}(0,\cdot)=0$ for all $n$.

\end{definition}

\begin{remark} \label{Hvaluedconti}
(i) Obviously, $\bH_{q,p}^{\alpha,\psi,\gamma+2}(T)$ is a Banach space.

(ii) By following the argument in \cite[Remark 3]{mikulevivcius2017p}, we can show that the embedding $H_p^{2n} \subset H_p^{\psi,2n}$ is continuous for any $n\in\bN$. This and Lemma \ref{H_p^phi,gamma space} (iv) imply that $\|u\|_{H^{\psi,2}_{p}} \leq C \|u\|_{H^{2}_{p}}$. Continuing the above argument, we obtain the desired result.
 \end{remark}

\begin{lemma} \label{basicproperty}
Let $\alpha\in(0,1)$, $1<p,q<\infty$, $\gamma\in\R$, and $T<\infty$.

(i) The space $\bH_{q,p,0}^{\alpha,\psi,\gamma+2}(T)$ is a closed subspace of $\bH_{q,p}^{\alpha,\psi,\gamma+2}(T)$.

(ii) $C_c^\infty(\bR^{d+1}_+)$ is dense in $\bH_{q,p,0}^{\alpha,\psi,\gamma+2}(T)$.

(iii) For any $\gamma,\nu\in\bR$, $(1-\mathcal{L})^{\nu/2}:\bH_{q,p}^{\alpha,\psi,\gamma+2}(T)\to\bH_{q,p}^{\alpha,\psi,\gamma-\nu+2}(T)$ is an isometry and for any $u\in \mathbb{H}^{\alpha,\psi,\gamma+2}_{q,p}(T)$, we have
$$
\partial^{\alpha}_{t}(1-\mathcal{L})^{\nu/2}u = (1-\mathcal{L})^{\nu/2}\partial^{\alpha}_{t}u. 
$$
\end{lemma}
\begin{proof}
If $\mathcal{L}=-\phi(-\Delta)$ for a Bernstein function $\phi$, then the lemma is a consequence of \cite[Lemma 2.7]{kim2020nonlocal}. Here we emphasize that the proof of \cite[Lemma 2.7]{kim2020nonlocal} is only based on Lemma \ref{H_p^phi,gamma space} and Remark \ref{Hvaluedconti} (ii). Hence, by following the proof, we can obtain all the claims.
\end{proof}

\begin{definition}\label{d:sv}
A function $f:(a,\infty)\to [0,\infty)$, for some $a>0$, is called slowly varying at infinity if for each $\lambda>0$
\begin{align*}
\lim_{x\to\infty}\frac{f(\lambda x)}{f(x)}=1.
\end{align*}
\end{definition}

Now, we introduce our assumption on the L\'evy density $j_{d}$.

\begin{assumption}\label{ass bernstein-2}
Let $\ell:(0,\infty) \to (0,\infty)$ be a continuous function. We assume that the jumping kernel $j_{d}:(0,\infty)\to(0,\infty)$ is differentiable and satisfies

\hspace{1mm} (i) The function $-\frac{1}{r}\left(\frac{d}{dr}j_{d}\right)(r)$ is decreasing;

\hspace{1mm} (ii) There exist constants $\kappa_{1},\kappa_{2}>0$ such that for any $m=0,1$ and $r>0$,
\begin{equation}\label{e:H}
\kappa_{1} r^{-d-m} \ell(r^{-1}) \leq (-1)^{m}\frac{d^m}{dr^m} j_{d}(r) \leq \kappa_{2} r^{-d-m}\ell(r^{-1}).
\end{equation}
\end{assumption}

The function $\ell$ gives the intensity of jumps of corresponding L\'evy process. Since $j_d(r)\asymp r^{-d}\ell(r^{-1})$ for $r>0$, the behavior of $\ell$ at infinity (resp. $0$) shows the behavior of small (resp. large) jumps.  We impose the following condition on $\ell$ in Assumption \ref{ass bernstein-2}.

\begin{assumption}\label{ass bernstein}
We assume that $\ell$ in Assumption \ref{ass bernstein-2} is independent of $d$ and  satisfies
\begin{align}
&C_{1}\left(\frac{R}{r}\right)^{\delta_{1}} \leq \frac{\ell(R)}{\ell(r)} \leq C_{2} \left( \frac{R}{r} \right)^{\delta_{2}} \quad \text{for}\;\; \,1\leq r\leq R<\infty;\label{H:s}
\end{align}
with constants $C_{1},C_{2}>0$ and  $0 \leq \delta_{1}\leq \delta_{2}<2$. We also assume that there exists $\delta_{3}>0$ such that
\begin{align}\label{H:l}
C_{1} \left( \frac{R}{r}  \right)^{\delta_{3}} \leq \frac{h(r)}{h(R)} \quad \text{for} \;\; \, 1\leq r\leq R<\infty,
\end{align} 
where 
\begin{equation}
\begin{gathered}\label{eqn 11.02.14:53}
h(r):= K(r) + L(r) \quad \text{for} \quad r>0,  
\\
K(r):=r^{-2}\int_{0}^{r}s\ell(s^{-1})ds, \quad L(r):=\int_{r}^{\infty} s^{-1}\ell(s^{-1}) ds \quad \text{for} \quad r>0.  
\end{gathered}
\end{equation}
\end{assumption}

In the rest of the article, we use a vector notation {\boldmath$\delta$} $=(\delta_{1},\delta_{2},\delta_{3})\in\bR^{3}$  instead of listing $\delta_{i}$ ($i=1,2,3$) for notational convenience.

\begin{remark}\label{r:assum}
(i) By \eqref{H:s} we obtain
\begin{align}\label{ll1}
\liminf_{r\to\infty}\ell(r)>0,
\end{align}
which is essential to our approach. If $\delta_1<0$,  \eqref{ll1} may not hold. For example, if $\ell(r)=\frac{1}{\log(1+r)}$, then $\delta_1$ should be negative. It is easy to see that \eqref{ll1} does not hold in this case. Thus, we only consider the case that $\delta_1\ge0$. On the other hand, \eqref{e:H} (with $m=0$) implies that $\ell$ bounded near zero since $j_{d}(|y|)dy$ is a L\'evy measure.

(ii)   Since the constant $\delta_1$ in \eqref{H:s} can be $0$, the function $\ell$ can be a slowly varying function at infinity. We can see  the relation between the characteristic exponent and the jump kernel of isotropic unimodal L\'evy process with low intensity of small jumps in  \cite{KM12,GRT19}.

(iii) Note that $K$, $L$ and $h$ are independent of the dimension $d$.
Clearly, $L$ is strictly decreasing in $r$.  Under \eqref{e:H}, for any $r>0$, we have
\begin{equation*}
\begin{gathered}
c(d)\kappa_{2}^{-1} r^{-2}\int_{|y|\leq r}|y|^{2} j_{d}(|y|)dy \leq K(r) \leq c(d) \kappa_{1}^{-1}r^{-2}\int_{|y|\leq r}|y|^{2} j_{d}(|y|)dy, 
\\
c(d) \kappa_{2}^{-1} \int_{|y|\geq r} j_{d}(|y|)dy \leq L(r) \leq c(d) \kappa_{1}^{-1} \int_{|y|\geq r} j_{d}(|y|)dy, 
\\
c(d) \kappa_{2}^{-1}r^{-2}\int_{\bR^{d}} \left( r^{2} \wedge |y|^{2} \right) j_{d}(|y|)dy \leq h(r) \leq c(d) \kappa_{1}^{-1}r^{-2}\int_{\bR^{d}} \left( r^{2} \wedge |y|^{2} \right) j_{d}(|y|)dy.
\end{gathered}
\end{equation*} 
Since $h'(r)=K'(r)+L'(r) = -2r^{-1}K(r)<0$, $h$ is strictly decreasing in $r$. Thus, the inverse function $h^{-1}$ of $h$ is well-defined and \eqref{H:l} makes sense. The function $\psi$ and $h$ are related as
\begin{equation}\label{eqn 05.27.15:40}
C_{0}h(r)\le \psi(r^{-1})\le 2h(r),\quad\text{for}\;\; r>0
\end{equation}
where the constant $C_{0}$ (we may take $C_{0}<1$) only depends on $d$ 
(see inequalities (6) and (7) in \cite{BGR14a}). Moreover, for any $0<c<1$, 
\begin{equation}\label{eqn 05.12:15:07}
h(cr) \leq c(d)\kappa_{1}^{-1} c^{-2} r^{-2} \int_{\bR^{d}} \left( c^{2} r^{2} \wedge |y|^{2} \right) j_{d}(|y|) dy  \leq \kappa_{2}\kappa_{1}^{-1} c^{-2} h(r).
\end{equation}

(iv) The condition \eqref{H:l} is about the behavior of $\ell$ near zero. Similar to \eqref{H:s}, we could give weak scaling conditions to $\ell$ near zero instead of \eqref{H:l}. However, to cover general decay rates for jump kernel, we choose to give condition \eqref{H:l}. Indeed, the behavior of $h$ near infinity is determined by the behavior of $\ell$ near zero only. To see this,  let $r>1$. Then we have
\begin{align}\label{eqn 10.28.17:51}
K(r) 
& \asymp  r^{-2}\int_{|y|\leq 1} |y|^{2} j_{d}(|y|)dy + r^{-2}\int_{1}^{r} s\ell(s^{-1})ds  \nonumber
\\
& \asymp  r^{-2} +r^{-2}\int_{1}^{r} s\ell(s^{-1})ds.
\end{align}
Thus, for $r>1$
$$h(r)=K(r)+L(r)\asymp r^{-2}+r^{-2}\int_{1}^{r} s\ell(s^{-1})ds+\int^{\infty}_{r}s^{-1}\ell(s^{-1})ds.$$
This shows that the behavior of $h$ near infinity depends only on the behavior of $\ell$ near zero. The following are some conditions on $\ell$ that give \eqref{H:l}: for $0<r\le R<1$
\begin{itemize}
\item $\ell(R)/\ell(r)\ge c(R/r)^{\delta_3'}$ with $\delta_3'>2$;
\item $c_1(R/r)^{\delta_3'}\le \ell(R)/\ell(r)\le c_2(R/r)^{\delta_4'}$ with $0<\delta_3'\le \delta_4'<2$;
\item $\ell(r)\asymp r^{2}$.
\end{itemize}

(v) Assumption \ref{ass bernstein-2} is used to obtain upper bounds for derivative of heat kernel. Using these upper bounds, we obtain upper bounds for derivative of our fundamental solution. If $X_{t}$ is a subordinate Brownian motion, then 
$$
j_{d+2}(r)=-\frac{1}{2\pi^{1/2}r}\frac{d}{dr}j_d(r) \quad \text{for}\; \, r>0
$$
holds due to \eqref{eqn 04.28.16:34}. Thus, Assumption \ref{ass bernstein-2} (i) holds for all subordinate Brownian motions. Suppose that for any $d\ge1$, the jump kernel $j_d$ of $d$-dimensional subordinate Brownian motion satisfies $j_d(r)\asymp r^{-d}\ell(r^{-1})$ for $r>0$ with a function $\ell$ independent of $d$. Then,  Assumption \ref{ass bernstein-2} (ii) also holds true.
\end{remark}

\begin{remark}\label{rmk 07.27.10:03}
Due to Assumption \ref{ass bernstein-2} we can apply \cite[Theorem 1.5]{kul2016gradient} to obtain $(d+2)$-dimensional isotropic unimodal L\'evy process $\tilde{X}_{t}$ with the same characteristic exponent $\psi(|\xi|)=\psi_X(|\xi|)$, whose transition density $p_{d+2}(t,x)=p_{d+2}(t,|x|)$ is radial,  radially decreasing in $x$ and satisfies
\begin{equation*}
p_{d+2}(t,r) = -\frac{1}{2\pi r} \frac{d}{dr}p_{d}(t,r) \quad\text{for}\;\; r>0.
\end{equation*}
This implies that for $t>0$ and $x\in\R^d$,
\begin{equation}\label{p'_est}
|D_x p_{d}(t,x)|  \leq 2\pi |x| p_{d+2}(t, |x|).
\end{equation}
By inspecting the proof of \cite[Theorem 1.5]{kul2016gradient}, we can also find that the jumping kernel $\tilde{j}_{d+2}$ of $\tilde{X}_{t}$ is given by
\begin{align}\label{j'_est}
\tilde{j}_{d+2}(r) = -\frac{1}{2\pi r} \frac{d}{dr} j_{d}(r) \quad\text{for}\;\; r>0.
\end{align}
By Assumption \ref{ass bernstein-2}, $\tilde{j}_{d+2}$ satisfies \eqref{e:H} for $m=0$ with $d+2$ in place of $d$. Thus, we can obtain the upper bounds for $|D_xp_{d}(t,x)|$ by using upper bounds of $p_{d+2}$.
\end{remark}

\medskip

\begin{remark}\label{rmk 05.28:15:16}\label{r:relations}
(i) Since $j_{d}(r)$ is decreasing in $r$, we have that for $r>0$
\begin{equation*}
\begin{aligned}
\ell(r^{-1}) \leq \kappa_{1}^{-1} r^{d} j_{d}(r) \leq (d+2) \kappa_{1}^{-1}  r^{-2}  \int_{0}^{r} s^{d+1} j_{d}(s) ds 
 \leq C K(r) \leq Ch(r),
\end{aligned}
\end{equation*}
where the constant $C$ depends only on $\kappa_{1},d$. 

(ii) From \eqref{H:l} and the fact that $h$ is strictly decreasing, we have the following: for any $0<A<\infty$, there exists a constant $c=c(A)$ such that
\begin{equation}\label{eqn 08.30.18:01}
c(A) \left(\frac{r^{-1}}{R^{-1}}\right)^{\delta_3} \leq \frac{h(r)}{h(R)}, \quad\text{for}\;\; A<r<R<\infty.
\end{equation}
Therefore, we have (put $r=h^{-1}(R)$, $R=h^{-1}(r)$)
\begin{equation}\label{eqn 05.27.16:40}
\frac{h^{-1}(r)}{h^{-1}(R)} \leq c(A) \left( \frac{r}{R} \right)^{1/\delta_3}, \quad\text{for}\;\; 0<r<R<h(A).
\end{equation}
Also, from \eqref{eqn 05.12:15:07}, we can obtain
\begin{equation}\label{eqn 05.30.17:20}
\frac{R}{r} \leq \left( \frac{h^{-1}(r)}{h^{-1}(R)} \right)^{2}, \quad\text{for}\;\; 0<r<R<\infty.
\end{equation}
\end{remark}

Now we introduce second assumption on $\ell$. Depending on whether $\ell$ is bounded or not, we have two different type of heat kernel upper bounds. Recall the due to Remark \ref{r:relations} (i), $\ell(r^{-1}) \leq C h(r)$ for all $r>0$.
 
\begin{assumption}\label{ell_con}
The function $\ell$ in Assumption \ref{ass bernstein} satisfies
\begin{itemize}
\item[(i)] either $\limsup_{r\to\infty}\ell(r)<\infty$;

\item[(ii)] or $ \limsup_{r\to\infty}\ell(r)=\infty$ and  $\ell(r)\asymp\sup_{s\le r}\ell(s)$.

\vspace{2mm}

\noindent If $\ell$ satisfies (ii), then we further assume that 

(ii)--(1) either
\begin{equation}\label{eqn 08.12.14:22}
\limsup_{r\to0}\frac{h(r)}{\ell(r^{-1})}<\infty.
\end{equation}

(ii)--(2) or 
\begin{equation*}\label{a22}
\limsup_{r\to0}\frac{h(r)}{\ell(r^{-1})}=\infty
\end{equation*}
and for any $a>0$ there is a constant $C(a)>0$ such that
\begin{equation}\label{eqn 08.12.17:59}
\sup_{0<r<1} h(r)\exp{\left( -a\frac{h(r)}{\ell(r^{-1})}  \right)} \leq C(a).
\end{equation}

\end{itemize}
\end{assumption}

\begin{remark}\label{rmk 10.24.17:53}
(i) If $0<\delta_{1} \leq \delta_{2}<2$ and \eqref{H:s} holds for all $0<r<R<\infty$, then we see that \eqref{eqn 08.12.14:22} holds.

(ii) Since $h$ is decreasing, \eqref{eqn 08.12.17:59} is equivalent to 
$\sup_{r>0} h(r)e^{-ah(r)/\ell(r^{-1})} \leq C(a)$. 

(iii) By \cite[Lemma 2.1]{cho21heat}, there exists $c>0$ such that
\begin{equation}\label{eqn 08.31.14:07}
L(r)\le h(r)\le cL(r)\quad\text{for}\;\;0<r\le1.
\end{equation}
By the change of variable, we see that $L(r)=\int^{r^{-1}}_{0}s^{-1}\ell(s)ds$ for $r>0$.
From Remark \ref{r:assum} (iii) and the fact that $j_{d}(|y|)dy$ is a L\'evy measure, we see that $L(1)=\int^{1}_{0}s^{-1}\ell(s)ds<\infty$. Since $L$ is strictly decreasing, 
\begin{align*}
L(r^{-1}) &
 = L(1) + \int_{1}^{r} s^{-1} \ell(s) ds 
\leq \frac{L(1)}{L(1/2)} L(r^{-1}) + \int_{1}^{r}s^{-1} \ell(s) ds \quad \text{for}\;\;r>2.
\end{align*}
Thus, $L(r^{-1})\asymp\int^{r}_{1}s^{-1}\ell(s)ds$ for $r>2$. Using this with \eqref{eqn 08.31.14:07} and (ii), it follows that \eqref{eqn 08.12.17:59} is equivalent to the following: for any $a>0$, there exists $C(a)>0$ such that
\begin{equation}\label{con_ii}
\sup_{r>1} \int^{r}_{1}\frac{\ell(s)}{s}ds\cdot\exp{\left( -\frac{a}{\ell(r)}\int^{r}_{1}\frac{\ell(s)}{s}ds  \right)} \leq C(a).
\end{equation}
It is known that for a slowly varying function $\ell:(b,\infty)\to(0,\infty)$, it holds that $\frac{1}{\ell(r)}\int^{r}_{b}\frac{\ell(s)}{s}ds\to\infty$ as $r\to\infty$ (see \cite[Proposition 1.5.9a]{BGT}). 
\end{remark}

Let $\cG$ be the set of functions $\ell$ satisfying \eqref{con_ii}. We can check that $(\log{(1+r)})^{k} \in \cG$ for any $k\in\bN$.  See Lemma \ref{lem 09.01.17:25} for more properties of $\cG$.

\vspace{2mm}

The following theorem is main result of this article. Note that to prove the parabolic ($\alpha=1$) counterpart of our main result, it seems that we need more differentiability to $j_{d}$. See Remark \ref{rmk 10.10.17:10} for detail.

\begin{theorem} \label{main theorem2}
Let $\alpha\in(0,1)$, $p,q\in(1,\infty)$, $\gamma \in \bR$, and $T\in(0,\infty)$. Suppose   Assumption \ref{ass bernstein-2},  Assumption \ref{ass bernstein} and Assumption \ref{ell_con} hold. Then for any  $f\in \bH_{q,p}^{\psi,\gamma}(T)$, the equation
\begin{equation}\label{mainequation2}
\partial_t^\alpha u = \mathcal{L}u + f,\quad t>0\,; \quad u(0,\cdot)=0
\end{equation}
has a unique solution $u$ in the class $\bH_{q,p,0}^{\alpha,\psi,\gamma+2}(T)$, and for the solution $u$ it holds that
\begin{equation} \label{mainestimate2}
\|u\|_{\bH_{q,p}^{\alpha,\psi,\gamma+2}(T)}\leq C \|f\|_{\bH_{q,p}^{\psi,\gamma}(T)},
\end{equation}
where $C>0$ depends only on $\alpha,d,\kappa_{1},\kappa_{2},p,q,\gamma,\ell,C_{0},C_{1},C_{2},T$ and {\boldmath$\delta$}. Furthermore, we have
\begin{equation}\label{eqn 05.27.14:12}
\| \mathcal{L} u\|_{\mathbb{H}^{\psi,\gamma}_{q,p}(T)} \leq C \| f\|_{\mathbb{H}^{\psi,\gamma}_{q,p}(T)},
\end{equation}
where $C>0$ depends only on $\alpha,d,\kappa_{1},\kappa_{2},p,q,\gamma,\ell,C_{0},C_{1},C_{2}$ and {\boldmath$\delta$}.
\end{theorem}

\begin{remark}\label{rmk 11.01.17:18} 
In Section \ref{sec 11.01.14:33}, we will discuss L\'evy processes such that $\ell$ depends on $d$. One of such examples is $\mathcal{L}=-\log{(1-\Delta)}$.
\end{remark}

\begin{remark}
In \cite{kim2020nonlocal}, \eqref{eqn 09.25.19:08} with non-zero initial data is considered. By using our result and following the approach in \cite{kim2020nonlocal}, we obtain the following: Let $\psi(r) = \phi_{k}(r^{2})$, where $\phi_k(r)=\log(1+r^{\beta/k})^k$ with $\beta\in(0,1)$ and $k\in\bN$. Then, the solution $u$ to \eqref{eqn 09.25.19:08} satisfies that 
\begin{equation*}\label{eqn 11.07.19:30}
\|u\|_{\bH_{q,p}^{\alpha,\psi,\gamma+2}(T)}\leq C  \left( \|u_{0}\|_{B^{\psi,\gamma+2+2/kq'-2/\alpha q}_{p,q}}  +  \|f\|_{\bH_{q,p}^{\psi,\gamma}(T)}  \right).
\end{equation*}
Compared to the result in \cite{kim2020nonlocal}, we need the additional regularity $2/kq'$ for initial data $u_0$ since $\phi_k$ does not satisfy \eqref{phi_lsc} (see \cite[Definition 2.3]{kim2020nonlocal} for $B^{\psi,s}_{p,q}$). It is nontrivial to remove this additional regularity for $u_0$.
\end{remark}

\begin{example}\label{exm 10.26.16:17}

The following are examples of $\ell$, $h$ and $j_d$ satisfying Assumption \ref{ass bernstein-2}, Assumption \ref{ass bernstein} and  Assumption \ref{ell_con}.

(i) Let $\beta\in(0,2)$, $\ell(r)=r^{\beta}$, and let $j_{d}(r) = r^{-d}\ell(r^{-1})$. Then $\ell$ satisfies \eqref{H:s} with $\delta_{1}=\delta_{2}=\beta$ and $K(r)\asymp L(r)\asymp h(r) \asymp r^{-\beta}$. Hence, we can easily check that Assumption \ref{ass bernstein-2} and Assumption \ref{ass bernstein} hold. Also, $\ell$ satisfies Assumption \ref{ell_con}(ii)-(1). Recall that the fractional Laplacian $-(-\Delta)^{\beta/2}$ is an example of operators covered by $\ell$ with explicit form. 
%Note that Assumption \ref{asm 11.07.16:00} holds with $\varepsilon_0=0$ by Remark \ref{r:u0}.

(ii) Let $\beta\in(0,2)$, $\ell(r) = r^{\beta}/(1+r^{\beta})$, and let $j_{d}(r)=r^{-d}\ell(r^{-1})$. Then $\ell$ satisfies \eqref{H:s} with $\delta_{1}=\delta_{2}=0$. Also, by a direct computation, we see that
\begin{align*}
-\frac{d}{dr}j_{d}(r) = d r^{-d-1}\frac{r^{-\beta}}{1+r^{-\beta}} + \beta r^{-d-1}\frac{r^{-\beta}}{(1+r^{-\beta})^{2}} \asymp r^{-d-1}\ell(r^{-1}) \quad \text{for}\quad  r>0,
\end{align*}
and thus $-\frac{1}{r}\left(\frac{d}{dr}j_{d} \right) (r)$ is a decreasing function (recall $\beta<2$). Since $\ell(s^{-1}) \asymp (1\wedge s^{-\beta})$, it is easy to see that
$L(r)  \asymp r^{-\beta}$ for $r\geq1$, and $K(r) \asymp r^{-2}+r^{-\beta} \asymp r^{-\beta}$ by \eqref{eqn 10.28.17:51}.  Thus $h(r) \asymp r^{-\beta}$ for $r\geq1$ and \eqref{H:l} follows with $\delta_{3}=\beta$. Hence,  $\ell$, $h$ and $j_{d}$ satisfy Assumption \ref{ass bernstein-2} and Assumption \ref{ass bernstein}. Note that, $\ell$ obviously satisfies Assumption \ref{ell_con}(i) since $\lim_{r\to\infty}\ell(r)=1<\infty$. 
%Moreover, Assumption \ref{asm 11.07.16:00} holds with $\varepsilon_0=1$ by Remark \ref{r:u0}. 

An operator $\mathcal{L}$ covered by the above function $\ell$ with explicit form is $\mathcal{L}=-\log{(1+(-\Delta)^{\beta/2})}$. In this case, the characteristic exponent is given as $\psi(r) =: \phi(r^{2}) =: \log{(1+r^{\beta})}$.  Then, it is known that $\phi(r) = \log{(1+r^{\beta/2})}$ is a Bernstein function and for each $d\geq1$, the corresponding $d$-dimensional process $X_{t}$ is a subordinate Brownian motion, with jumping kernel $J_{d}$ satisfying (see e.g. \cite{KM12,KSV14})
$$
J_{d}(r) \asymp r^{-d}(1\wedge r^{-\beta}) \asymp r^{-d} \frac{r^{-\beta}}{1+r^{-\beta}} \quad \text{for}\;\, r>0. 
$$
Hence, we deduce that $J_{d}$ satisfies  Assumption \ref{ass bernstein-2} with $\ell(r)=r^{\beta}/(1+r^{\beta})$.

(iii) Let $\beta>0$, $\ell(r) =   \log{(1+r^{\beta})}$, and let $j_{d}(r)= r^{-d}\ell(r^{-1})$. Then $\ell$ satisfies \eqref{H:s} with $\delta_{1}=\delta_{2}=0$.  Using the argument in (i) we can check that for $r\geq 1$
\begin{align*}
h(r) \asymp \begin{cases} r^{-\beta} & \text{if} \quad \beta\in(0,2)
\\
r^{-2}\log{r} & \text{if} \quad \beta=2
\\
r^{-2} & \text{if} \quad \beta>2.
\end{cases}
\end{align*}
Hence, Assumption \ref{ass bernstein} follows. Observe that 
\begin{align*}
-\frac{1}{r}\frac{d}{dr}j_d(r)
=dr^{-d-2}\log{(1+r^{-\beta})}+\beta r^{-d-2}\frac{1}{1+r^{\beta}}.
\end{align*}
Thus $-\frac{1}{r} \frac{d}{dr}j_{d} (r)$ is a decreasing function. Moreover, using 
\begin{align}\label{eqn 10.26.15:30}
(1+r)\log(1+r^{-1})\ge\frac12 \quad \text{for}\;\; r>0,
\end{align}
we can check that $j_{d}$ satisfies Assumption \ref{ass bernstein-2}. Also, using Lemma \ref{lem 09.01.17:25} (ii), (iv) and Remark \ref{rmk 10.24.17:53} (iii), we see that $\ell$ satisfies Assumption \ref{ell_con}(ii)-(2).

The operator $\mathcal{L}=-(\log{(1+(-\Delta)^{\beta/4})})^{2}$ with $\beta\in(0,2)$ can be covered by $\ell(r) = \log{(1+r^{\beta})}$. In this case, $\psi(r) =: \phi(r^{2}) =: (\log{(1+r^{\beta/2})})^{2}$, and by using theory of Bernstein functions in \cite[Chapter 7]{SSV12}, we see that $\phi$ satisfies conditions (A-2) and (A-3) in \cite{KM12}. Hence, by \cite[Proposition 4.1]{KM12}, we can check that for each $d\geq1$, jumping kernel $J_{d}$ of corresponding $d$-dimensional subordinate Brownian motion satisfies
\begin{align*}%\label{eqn 10.25.17:44}
J_{d}(r) \asymp r^{-d-2}\phi'(r^{-2}) \asymp r^{-d}\log{(1+r^{-\beta})} &\quad \text{for} \;\, r\leq 1.
\end{align*} 
Since $\phi(r)$ satisfies condition (H2) in \cite{KSV14}, we use \cite[Lemma 3.3(a)]{KSV14} to obtain that for each $d\geq 1$ we have
\begin{align*}
J_{d}(r) \asymp r^{-d}\phi(r^{-2}) = r^{-d}(\log{(1+r^{-\beta/2})})^{2}
 \asymp r^{-d}\log{(1+r^{-\beta})} \quad \text{for} \;\, r\geq 1.
\end{align*}
Using these estimates for $J_d$, for each $d\geq1$, we have $J_{d}(r) \asymp r^{-d}\ell(r^{-1})$ for all $r>0$ and hence it satisfies Assumption \ref{ass bernstein}, Assumption \ref{ass bernstein-2} and Assumption \ref{ell_con}.

(iv) For $b\in(0, 1/2)$, let $\ell(r):=\exp\{(\log(1+r))^{b}\}-1$. Then, by \cite[Theorem 1.3.1]{BGT}(put $c(x)=1$ and $\varepsilon(x)=(\log{x})^{b-1}$ therein), we can check that $\ell$ is a slowly varying function at infinity. Moreover, $\ell(r)\asymp r^b$ for $r<1$. Thus, we see that $\ell$ and $h$ satisfy \eqref{H:s} and \eqref{H:l} respectively. Let $j_d(r):=r^{-d}\ell(r^{-1})$. Then, using \eqref{eqn 10.26.15:30} and the same argument as in (ii), it is easy to check that $j_d(r)$ satisfies Assumption \ref{ass bernstein-2} with $\ell$. Finally, due to Lemma \ref{lem 09.01.17:25} (v), $\ell$ satisfies Assumption \ref{ell_con} (ii)-(2).

\vspace{2mm}
The above examples can be summarized as follows.
\begin{center}
\renewcommand{\arraystretch}{1.4}%상하여백
\begin{tabular}[c]{| l | l | l |  c |}
	\hline
	$\ell(r)$    & $\psi(r) \asymp h(r^{-1})$ & Assumption \ref{ell_con}&  $(\delta_1, \delta_2, \delta_3)$\\
	\hline\hline
	$r^{\beta}$, $\beta\in(0,2)$     &    $\asymp r^{\beta}$  &    (ii)--(1)     &  $(\beta,\beta,\beta)$\\
	$r^{\beta}/(1+r^{\beta})$, $\beta\in(0,2)$   &  $\asymp \log (1+r^{\gamma})$ & (i)  & $(0,0,\beta)$ \\
	$\log (1+r^{\beta})$, $\beta\in(0,2)$   &  $\asymp (\log(1+r^{\beta/2}))^2$ & (ii)--(2)  & $(0,1,\beta)$\\
	\hline
\end{tabular}
\end{center}

\end{example}
\medskip

\mysection{Estimates of the heat kernels and their derivatives}
\label{sec 11.09.16:14}

In this section, we obtain  sharp bounds of  the heat kernel and its derivative for equation $\partial_t u=\mathcal{L}u$ under the Assumption \ref{ass bernstein-2}, Assumption \ref{ass bernstein} and Assumption \ref{ell_con}. For the rest of this section, we suppose that the Assumption \ref{ass bernstein-2} and Assumption \ref{ass bernstein} hold.

Let $p(t,x)=p_d(t,x)$ be the transition density of $X_t$. Then it is well-known that for any $t>0, x\in \bR^d$,
\begin{equation}\label{eqn 7.20.1}
p_d(t,x)=\frac{1}{(2\pi)^d} \int_{\bR^d} e^{i \xi \cdot x} e^{-t\psi(|\xi|)} \,d\xi.
\end{equation}
Since $X_{t}$ is isotropic, $p_{d}(t,x)$ is rotationally invariant in $x$ (i.e. $p_{d}(t,x)=p_{d}(t,|x|)$). We put $p_{d}(t,r):=p_{d}(t,x)$  if $r=|x|$ for notational convenience. Since $X_{t}$ is unimodal, $r\mapsto p(t,r)$ is a decreasing function. Moreover, $p(t,x)\le p(t,0)\in(0,\infty]$ for $t>0$ and $x\in\R^d$.

The first part of this section consists of estimates for the heat kernel $p_{d}$ to the equation 
\begin{equation}\label{mainequation}
\partial_t u = \mathcal{L}u.
\end{equation}
The following lemma gives off diagonal type upper bound for the heat kernel. The result holds for all isotropic unimodal L\'evy processes. 

\begin{proposition}[{\cite[Theorem 5.4]{GRT19}}]\label{p:uhk_jump}
For any $(t,x)\in (0,\infty)\times\bR^{d}$, we have
\begin{equation*}
p_{d}(t,x) \leq C t|x|^{-d}K(|x|),
\end{equation*}
where the constant $C>0$ depends only on $\kappa_{2},d$.
\end{proposition}

Next, we will give the sharp upper bound for heat kernel. Depending on whether $\ell$ is bounded or not, we have two different type of the upper bounds. If $\ell$ satisfies the condition in Assumption \ref{ell_con} (ii), we define 
$$\bar{\ell}(r):=\sup_{s\le r}\ell(s)\quad\text{for}\;\;r>0.$$
Then, $\bar{\ell}:(0,\infty)\to(0,\infty)$ is increasing continuous function and $\ell(r)\asymp \bar{\ell}(r)$ for all $r>0$. Also since  $\ell$ satisfies Assumption \ref{ass bernstein} and is bounded near zero (recall Remark \ref{r:assum} (i)), we can apply Lemma \ref{lem 07.23.16:47} to obtain a strictly increasing continuous function $\ell^{\ast}$ satisfying
\begin{equation*}
\ell^{\ast} \asymp \bar{\ell} \asymp \ell.
\end{equation*}
Let $\ell^{-1}$ be the inverse function of $\ell^{\ast}$.

For any $a>0$ and $r,t>0$, we define
\begin{align}\label{d:theta}
\theta(a,r,t)=\theta_a(r,t):=r\vee(\ell^{-1}(a/t))^{-1}.
\end{align}
The following heat kernel upper bounds are proved in  \cite[Proposition 2.9, Corollary 2.13]{cho21heat}.

\begin{proposition}\label{p:hku}
(i) Suppose $\ell$ satisfies the condition in Assumption \ref{ell_con} (i). For any $T>0$, there exist constants $c_{1}, b_{1}>0$ depending only on $\kappa_{1},\kappa_{2},d$ and $T$ such that
\begin{equation}\label{uhk_bdd}
p_{d}(t,x) \leq c_{1}t\frac{K(|x|)}{|x|^d}\exp{(-b_{1}th(|x|))}
\end{equation}
holds for all $(t,x)\in (0,T]\times (\bR^{d}\setminus\{0\})$. 

(ii) Suppose $\ell$ satisfies the condition in Assumption \ref{ell_con} (ii) and $\theta$ is given by \eqref{d:theta}.  
Then, there exist constants $a_0, b_2>0$ such that the following holds: for any $T>0$, there exist $c_{2}>0$ depending only on $\kappa_{1},\kappa_{2},d$ and $T$ such that
\begin{equation}\label{uhk_unbdd}
p_{d}(t,x) \leq c_{2}t\frac{K(\theta_{a_0}(|x|,t))}{[\theta_{a_0}(|x|,t)]^d}\exp{(-b_{2}th(\theta_{a_0}(|x|,t)))}
\end{equation}
holds for all $(t,x)\in (0,T]\times \bR^{d}$. 
\end{proposition}
\begin{proof}
If $\ell$ satisfies Assumption \ref{ell_con} (i), then \eqref{uhk_bdd} is a direct consequence of \cite[Corollary 2.13]{cho21heat}. If we define 
$$
\theta'(a,r,t)=\theta'_{a}(r,t)  :=r\vee(\bar{\ell}^{-1}(a/t))^{-1},
$$ 
where $\bar{\ell}^{-1}$ is the right continuous inverse  of $\bar{\ell}$, then \eqref{uhk_unbdd} holds with $\theta'$ in place of $\theta$ due to \cite[Corollary 2.9]{cho21heat}. Since $\ell^{\ast}\asymp \bar{\ell}$, by inspecting the proof of \cite[Corollary 2.9]{cho21heat}, we can check $\theta'$ in \cite[Corollary 2.9]{cho21heat} can be replaced by $\theta_{a}$. The lemma is proved. 
\end{proof}

We will only use the upper bound of $p_{d}(t,x)$ in \eqref{uhk_bdd} for the proofs of our results. 
The following lemma will be used several times.

\begin{lemma}\label{rmk 05.07.16:12}
Let $a, b>0$ and $d\in\bN$.\\
(i) There exists $C=C(b,d)>0$ such that for all $t>0$
\begin{equation*}
\begin{aligned}
\int_{\bR^{d}} t\frac{K(|x|)}{|x|^{d}}\exp{(-bth(|x|))}\,dx \leq C.
\end{aligned}
\end{equation*}
(ii) Suppose $\ell$ satisfies the condition in Assumption \ref{ell_con} (ii) and $\theta$ is given by \eqref{d:theta}. There exists $C=C(a,b,d)>0$ such that for all $t>0$
\begin{equation*}
\begin{aligned}
\int_{\bR^{d}} t\frac{K(\theta_{a}(|x|,t))}{[\theta_{a}(|x|,t)]^d}\exp{(-b th(\theta_{a}(|x|,t)))}\,dx \leq C.
\end{aligned}
\end{equation*}
\end{lemma}
\begin{proof}
(i) It is easy to see that
\begin{equation}\label{eqn 05.14.15:08}
\begin{aligned}
&\int_{\bR^{d}} t|x|^{-d}K(|x|) e^{-bth(|x|)} dx = C(d) \int_{0}^{\infty} t\rho^{-1} K(\rho) e^{-bth(\rho)} d\rho
\\
&=C \int_{0}^{h^{-1}(t^{-1})}  t\rho^{-1} K(\rho) e^{-bth(\rho)} d\rho +C \int_{h^{-1}(t^{-1})}^{\infty}  t\rho^{-1} K(\rho) e^{-bth(\rho)} d\rho.
\end{aligned}
\end{equation}
Note that for any $a>0$,
\begin{equation}\label{eqn 05.11.17:32}
\begin{aligned}
\int_{a}^{\infty}\rho^{-1}K(\rho) d\rho &= \int_{a}^{\infty}\rho^{-3}\int_{0}^{\rho}s \ell(s^{-1})ds d\rho
\\
& = \int_{0}^{a}\int_{a}^{\infty}\rho^{-3} d\rho s \ell(s^{-1}) ds + \int_{a}^{\infty}\int_{s}^{\infty}\rho^{-3}d\rho s \ell(s^{-1}) ds
\\
& = \frac{a^{-2}}{2} \int_{0}^{a} s \ell(s^{-1}) ds + \frac{1}{2} \int_{a}^{\infty}s^{-1} \ell(s^{-1})ds
\\
&=\frac{1}{2}\left( K(a)+L(a) \right) \leq h(a).
\end{aligned}
\end{equation}
Also, by the the change of variable $h(\rho)\to \rho$ and the fact that $h'(\rho)=-2\rho^{-1}K(\rho)$, we have
\begin{equation}\label{eqn 05.14.15:03}
\begin{aligned}
\int_{0}^{a} t \rho^{-1} K(\rho) e^{-bth(\rho)} d\rho \leq \frac{1}{2} t \int_{h(a)}^{\infty} e^{-b t \rho} d\rho = C e^{-bth(a)}.
\end{aligned}
\end{equation}
Applying \eqref{eqn 05.11.17:32} and \eqref{eqn 05.14.15:03} to \eqref{eqn 05.14.15:08} with $a=h^{-1}(t^{-1})$, we obtain the first assertion.

(ii) By the definition of $\theta_a(r,t)$, we see that
\begin{align*}
&\int_{\bR^{d}} t\frac{K(\theta_{a}(|x|,t))}{[\theta_{a}(|x|,t)]^d}\exp{\big(-b th(\theta_{a}(|x|,t))\big)}\,dx\\
&=\int_{|x|\ge (\ell^{-1}(a/t))^{-1}} t\frac{K(|x|)}{|x|^d}\exp{\big(-b th(|x|)\big)}\,dx\\
&\quad+\int_{|x|< (\ell^{-1}(a/t))^{-1}} t\ell^{-1}(a/t)^{d}{K(\ell^{-1}(a/t)^{-1})}\exp{\big(-b th(\ell^{-1}(a/t)^{-1})\big)}\,dx\\
&\leq \int_{\bR^{d}} t\frac{K(|x|)}{|x|^d}\exp{\big(-b th(|x|)\big)}\,dx\\
&\quad+\int_{|x|< (\ell^{-1}(a/t))^{-1}} t\ell^{-1}(a/t)^{d}{K(\ell^{-1}(a/t)^{-1})}\exp{\big(-b th(\ell^{-1}(a/t)^{-1})\big)}\,dx.
\end{align*}
The first term on the right hand side is bounded above by a constant $C$ due to (i). The second term can be easily handled by using the relations $K\leq h$ and $se^{-s}\leq 1$. Thus, we obtain the desired result.

\end{proof}

The following lemmas give upper bound of $p_{d}(t,x)$ for sufficiently large $t>0$ without  Assumption \ref{ell_con}.

\begin{lemma}\label{l:hke_largetime}
There exist $t_{1}=t_{1}(d,\kappa_{1},\kappa_{2},\ell,C_{0},C_{1},C_{2},\text{{\boldmath$\delta$}})>0$ and $C>0$ depending only on $t_{1}$ such that for all $t\geq t_{1}$ and $x\in\bR^{d}$,
\begin{align*}
p_d(t,x)\le C\left( (h^{-1}(t^{-1}))^{-d}\wedge t\frac{K(|x|)}{|x|^d}\right).
\end{align*}
\end{lemma}

\begin{proof}
First, we show that there is $t_{0}>0$ such that $e^{-C_{0}th(|\cdot|^{-1})}\in L_{1}$ for $t\geq t_{0}$, where $C_{0}$ is a constant from \eqref{eqn 05.27.15:40}. By  \eqref{eqn 05.27.15:40} and \eqref{eqn 08.31.14:07}, we see that  $\psi(r^{-1})\asymp L(r)$ for $r<1$. Hence, if $|\xi|\geq 1$, then we have
$$
\int_{|\xi|^{-1}}^{1} s^{-1} \ell(s^{-1})ds \leq L(|\xi|^{-1}) \leq c_1 \psi(|\xi|).
$$
Thus, by \eqref{ll1} and the above inequality, there exists $C>0$ such that
$$C\log{|\xi|}\le \psi(|\xi|)\quad  \text{for}\; \; |\xi|>1.$$
Thus if $t>0$ satisfies $C t \geq 2d$, then we have 
$$
e^{-t\psi(|\xi|)} \leq e^{-Ct\log{|\xi|}} \leq |\xi|^{-2d} \quad  \text{for}\; \; |\xi|>1.
$$
This and \eqref{eqn 05.27.15:40} show that there is $t_{0}>0$ such that $e^{-C_{0}t h(|\cdot|^{-1})}\in L_{1}$ for $t\geq t_{0}$.

On the other hand, by \eqref{eqn 08.30.18:01} we obtain that (recall $C_{0}$ comes from \eqref{eqn 05.27.15:40})
\begin{equation*}
c_2 \left(\frac{|\xi|}{(h^{-1}(t^{-1}))^{-1}}\right)^{\delta_{3}} \leq C_{0}th(|\xi|^{-1}) \quad\text{for}\;\; (h^{-1}(t^{-1}))^{-1}\leq |\xi| \leq 1.
\end{equation*}
Take $t_1\ge 2t_0$ satisfying $1 \leq h^{-1}(t_1^{-1})$. Then, using the above observation, we have that for $t\ge t_1$ and $|\xi|\geq1$,
\begin{align*}
e^{-C_{0}th(|\xi|^{-1})}
\leq e^{-C_{0}t_{0}h(|\xi|^{-1})}e^{-C_{0}t h(1)/2}
\le e^{-C_{0}t_{0}h(|\xi|^{-1})}e^{-c_2(h^{-1}(t^{-1}))^{\delta_{3}}}.
\end{align*}
Using the above inequalities, \eqref{eqn 7.20.1}, and \eqref{eqn 05.27.15:40}, for $t\geq t_{1}$, we obtain
\begin{equation*}
\begin{aligned}
|p_{d}(t,x)| &\leq \int_{\bR^{d}} e^{-C_{0} t h(|\xi|^{-1})} d\xi 
\\
 &\leq  \int_{|\xi| \leq (h^{-1}(t^{-1}))^{-1}}  1 d\xi + \int_{(h^{-1}(t^{-1}))^{-1} \leq |\xi| \leq 1}  e^{-c_2|\xi|^{\delta_{3}}(h^{-1}(t^{-1}))^{\delta_{3}}}d\xi 
\\
&\ \quad + \int_{|\xi|\geq1} e^{-C_{0}t_{0}h(|\xi|^{-1})}e^{-c_2(h^{-1}(t^{-1}))^{\delta_{3}}} d\xi
\\
&\leq c_3 (h^{-1}(t^{-1}))^{-d} + (h^{-1}(t^{-1}))^{-d} \int_{1\leq |\xi|} e^{-c_{2}|\xi|^{\delta_{3}}} d\xi
\\
&\quad + c_3 (h^{-1}(t^{-1}))^{-d} \int_{|\xi|\geq1} e^{-C_{0}t_{0}h(|\xi|^{-1})} d\xi
\\
&\leq c_4 (h^{-1}(t^{-1}))^{-d},
\end{aligned}
\end{equation*}
where for the third and last inequalities, we used that there exists $c_5(c_2, \delta_{3})>0$ such that $e^{-c_{2}x^{\delta_{3}}}\le c_5x^{-d}$ for all $x\geq 1$ and integrability of $e^{-t_{0}h(|\cdot|^{-1})}$, respectively. This completes the proof.
\end{proof}

\begin{remark}\label{r:hku_largetime}
(i) In Lemma \ref{l:hke_largetime}, we use \eqref{ll1} to show that $e^{-C_{0}th(|\cdot|^{-1})}$ is integrable for $t\geq t_{0}$. If we assume that $\ell$ satisfies Assumption \ref{ell_con} (ii), then we obtain the same near diagonal upper bound without \eqref{ll1}. Indeed, if $\ell$ satisfies Assumption \ref{ell_con} (ii), then by Proposition \ref{p:hku} (ii), we see that for fixed $T>0$
\begin{align}\label{eqn 08.18.12:17}
p_d(t,0)=\int_{\R^d}e^{-t\psi(|\xi|)} d\xi\le C\ell^{-1}(a_0/t)^{d}e^{-crh(\ell^{-1}(a_0/t)^{-1})} \quad \text{for}\;\;t\le T.
\end{align}
This implies that 
\begin{align*}
\int_{\R^d}e^{-t\psi(|\xi|)} d\xi\le \int_{\R^d}e^{-T\psi(|\xi|)} d\xi\le C\ell^{-1}(a_0/T)^{d}e^{-crh(\ell^{-1}(a_0/T)^{-1})} \quad \text{for}\;\;t>T.
\end{align*}
Thus, $e^{-t\psi(|\cdot|)}$ is integrable for all $t>0$.

(ii) Suppose that the function $\ell$ satisfies Assumption \ref{ell_con} (ii)--(1). Then by using Lemma \ref{l:hke_largetime} and  \eqref{eqn 08.18.12:17}, we can check that
$$
|p_{d}(t,x)| \leq |p_{d}(t,0)| \leq C (h^{-1}(a_{0}/t))^{-d} \leq C (h^{-1}(t^{-1}))^{-d} \quad \text{for}\; t>0,
$$
where the last inequality follows from \eqref{eqn 05.30.17:20}. Hence, if $\ell$ satisfies Assumption \ref{ell_con} (ii)--(1), then
\begin{align}\label{hke_largetime2}
p_d(t,x)\le C\left( (h^{-1}(t^{-1}))^{-d}\wedge t\frac{K(|x|)}{|x|^d}\right).
\end{align}
holds for all $(t,x)\in(0,\infty)\times\bR^{d}$.
\end{remark}

\medskip

\mysection{Estimates of the fundamental solutions}\label{sec 4}
Let $Q_t$ be an increasing L\'evy process independent of $X_t$ having the Laplace transform
\begin{equation*}
\bE \exp(-\lambda Q_t)= \exp(-t\lambda^\alpha),
\end{equation*}
and let 
$$
R_t:=\inf\{s>0 : Q_s>t\}
$$ 
be the inverse process of the  subordinator $Q_{t}$. Denote $\varphi(t,r)$ denote the probability density function of $R_t$ and define 
\begin{eqnarray*}
q_{d}(t,x):=\int_0^\infty p_{d}(r,x)d_r \bP(R_t\leq r)
=\int_0^\infty p_{d}(r,x) \varphi(t,r) \,dr.  
\end{eqnarray*}
We will show in Lemma \ref{zero converge} that the function $q_{d}(t,x)$ becomes the fundamental solution to
\begin{equation}\label{eqn 06.15.16:05}
\partial^{\alpha}_{t}u = \mathcal{L} u \quad t>0, \quad u(0,\cdot)=u_{0}.
\end{equation}
In this section, we also provide a sharp estimation of $q_{d}$ and its derivatives. Note that we can find a result for subordinate Brownian motion in \cite[Lemma 5.1]{kim2020nonlocal} (see also \cite[Theorem 1.1]{chen2017time}).

For $\beta\in \bR$, denote
$$
 \varphi_{\alpha,\beta}(t,r):=D_t^{\beta-\alpha}\varphi(t,r):=(D^{\beta-\alpha}_t \varphi(\cdot,r))(t).
$$
It is known (see e.g. \cite[Lemma 3.7 (ii)]{kim2020nonlocal}) that
\begin{equation} \label{philarge}
|\varphi_{\alpha,\beta}(t,r)|\leq C t^{-\beta}e^{-c(rt^{-\alpha})^{1/(1-\alpha)}}
\end{equation}
for $rt^{-\alpha}\geq 1$, and 
\begin{equation}\label{betainteger}
|\varphi_{\alpha,\beta}(t,r)|\leq \left\{ \begin{array}{ll}  C rt^{-\alpha-\beta}~&\beta\in \bN \\ Ct^{-\beta}~&\beta \notin \bN\end{array} \right.
\end{equation}
for $rt^{-\alpha}\leq 1$, where the constants $c,C>0$ depending only on $\alpha,\beta$.

For  $(t,x)\in(0,\infty)\times \bR^{d}\setminus\{0\}$ define
\begin{equation*}
q^{\alpha,\beta}_{d}(t,x):=\int_0^\infty p_{d}(r,x)\varphi_{\alpha,\beta}(t,r) dr.
\end{equation*}
Note that since $\varphi_{\alpha,\beta}(t,r)$ is integrable in $r\in(0,\infty)$, we can check that $q^{\alpha,\beta}_{d}$ is well-defined.

\begin{lemma}
 \label{lem 7.21.1}
 Let $\alpha\in (0,1)$ and $\beta\in \bR$.

(i) For  any $(t,x)\in(0,\infty)\times \bR^{d}\setminus\{0\}$, 
\begin{equation} \label{eqn 07.28 20:25}
D_t^{\beta-\alpha}q_{d}(t,x)=q^{\alpha,\beta}_{d}(t,x).
\end{equation}

(ii) For any $t>0, \xi\in \bR^d$,
\begin{equation} \label{fourier of q}
\mathcal{F}\{q^{\alpha,\beta}_{d}(t,\cdot)\}(\xi) =t^{\alpha-\beta} E_{\alpha,1-\beta+\alpha}(-t^\alpha\psi(|\xi|)),
\end{equation}
where $E_{a,b}$ be the two-parameter Mittag-Leffler function defined as
\begin{equation*}
E_{a,b}(z)=\sum_{k=0}^\infty \frac{z^k}{\Gamma(a k+b)}, \quad\text{for}\;\;z\in\bC, a>0, b\in\bC,
\end{equation*}
with the convention $E_{a,1}(0)=1$.
\end{lemma}

\begin{proof}
(i) By Proposition \ref{p:uhk_jump}, we have 
\begin{equation*}
	\label{eqn 07.28.22:31}
 |p_{d}(r,x)| \leq C(d,\kappa_{2},x, K) r.
\end{equation*}
Hence, we can prove the assertion by following the proof of \cite[Lemma 3.7 (iii)]{kim2020nonlocal}.

(ii) We can directly prove this by following the proof of \cite[Lemma 3.7 (iv)]{kim2020nonlocal} with $\psi(|\xi|)$ in place of $\phi(|\xi|^{2})$ since we need no property of $\phi$.

The lemma is proved.
\end{proof}

For the rest of the paper, we assume that Assumption \ref{ass bernstein}, Assumption \ref{ass bernstein-2} and Assumption \ref{ell_con} hold.

\begin{lemma}\label{lem 05.27.18:14}
Let $\alpha\in(0,1)$, $\beta\in\bR$, and $m=0,1$. Then there exists a constant $C$ depending only on $\alpha,\beta,\kappa_{2},d$ and $m$ such that for any $t>0$ and $x\in \bR^{d}\setminus\{0\}$,
\begin{equation}\label{eqn 07.02.16:48}
|D^{m}_{x}q^{\alpha,\beta}_{d}(t,x) | \leq C t^{2\alpha-\beta} \frac{K(|x|)}{|x|^{d+m}}.
\end{equation}
\end{lemma}

\begin{proof}
Using \eqref{p'_est} and Proposition \ref{p:uhk_jump}, we have that for any $x\in \bR^{d}\setminus\{0\}$, $r>0$ and for any $y\in B_{\varepsilon}(x)$ for sufficiently small $\varepsilon>0$, 
\begin{align}\label{pd_diff}
|D^{\sigma}_{x}p_{d}(r,y)| 
\le C( \kappa_2, d)r\frac{K(|y|)}{|y|^{d+|\sigma|}}
\le C(\kappa_{2},d,x,\varepsilon,K) r \;\; \text{for}\; |\sigma| \leq 1.
\end{align}
Hence, by the dominated convergence theorem, we have
\begin{align}\label{eqn 09.03.19:06}
D^{m}_{x}q^{\alpha,\beta}_{d}(t,x) = \int_{0}^{\infty} D^{m}_{x}p_{d}(r,x) \varphi_{\alpha,\beta}(t,r)dr.
\end{align}
By \eqref{pd_diff}, \eqref{philarge} and \eqref{betainteger}
\begin{equation*}
\begin{aligned}
|D^{m}_{x}q^{\alpha,\beta}_{d}(t,x)| &\leq C\int_{0}^{t^{\alpha}} rK(|x|) |x|^{-d-m} t^{-\beta} dr
\\
& \quad + C \int_{t^{\alpha}}^{\infty} rK(|x|) |x|^{-d-m} t^{-\beta} e^{-c(rt^{-\alpha})^{\frac{1}{1-\alpha}}} dr
\\
& \leq C t^{2\alpha-\beta} \frac{K(|x|)}{|x|^{d+m}},
\end{aligned}
\end{equation*}
where for the second integral, we used the change of variable $rt^{-\alpha} \to r$. Thus, the lemma is proved.
\end{proof}

\begin{lemma}\label{l:qest_largetime}
Let $\alpha\in(0,1)$, $\beta\in\bR$, and $m=0,1$, and let $t_{1}>0$ be taken from Lemma \ref{l:hke_largetime}. Then for any $(t,x)\in(0,\infty)\times\bR^{d}$ satisfying $t^{\alpha}h(|x|) \geq 1$, we have
\begin{align*}
&  \int_{t_1}^{\infty} |D^{m}_{x}p_{d} (r,x)| |\varphi_{\alpha,\beta} (t,r)| dr  \leq C  \int_{(h(|x|))^{-1}}^{2t^{\alpha}} (h^{-1}(r^{-1}))^{-d-m} r t^{-\alpha-\beta} dr 
\end{align*}
for $\beta\in\bN$,
and
\begin{align*}
&  \int_{t_1}^{\infty} |D^{m}_{x}p_{d} (r,x)| |\varphi_{\alpha,\beta} (t,r)| dr  \leq C  \int_{(h(|x|))^{-1}}^{2t^{\alpha}} (h^{-1}(r^{-1}))^{-d-m}  t^{-\beta} dr 
\end{align*}
for $\beta\notin\bN$, where $C$ depends only on $\alpha,\beta,d,\kappa_{1},\kappa_{2},\text{\boldmath{$\delta$}},\ell,C_{0},C_{1},C_{2}$ and $m$. If $\ell$ satisfies Assumption \ref{ell_con} (ii)--(1), then we can put $t_{1}=0$ in the above inequalities.
\end{lemma}

\begin{proof}
Suppose that $\beta\in\bN$ and $h(|x|) \leq t_{1}^{-1}$.  Then, by Remark \ref{rmk 07.27.10:03}, Lemma \ref{l:hke_largetime},  \eqref{philarge} and \eqref{betainteger}, we have
\begin{equation*}
\begin{aligned}
\left|  \int_{t_1}^{\infty} |D^{m}_{x}p_{d} (r,x)| |\varphi_{\alpha,\beta} (t,r)| dr \right|
  \leq C\left( I_{1}+I_{2}+I_{3} \right),
\end{aligned}
\end{equation*}
where
\begin{align}
I_{1} &= \int_{t_{1}}^{(h(|x|))^{-1}} r \frac{K(|x|)}{|x|^{d+m}} r t^{-\alpha-\beta}dr,
\nn\\
I_{2} &=  \int_{(h(|x|))^{-1}}^{t^{\alpha}} |x|^{m} (h^{-1}(r^{-1}))^{-d-2m} r t^{-\alpha-\beta} dr,
\nn\\
I_{3} &=  \int_{t^{\alpha}}^{\infty} |x|^{m} (h^{-1}(r^{-1}))^{-d-2m}  t^{-\beta} e^{-c(rt^{-\alpha})^{\frac{1}{1-\alpha}}} dr, \label{eqn 05.28.12:51}
\end{align}
and the constant $C$ depends only on $\alpha,\beta,d,\kappa_{1},\kappa_{2},\text{\boldmath{$\delta$}},\ell,C_{0},m$. We first estimate $I_{1}$. It is easy to see that
\begin{equation*}
\begin{aligned}
I_{1} \leq (h(|x|))^{-2} |x|^{-d-m} t^{-\alpha-\beta} = C \int_{(h(|x|))^{-1}}^{2(h(|x|))^{-1}}  |x|^{-d-m} r t^{-\alpha-\beta}dr.
\end{aligned}
\end{equation*}
For $t_{1}<(h(|x|))^{-1}<r<2(h(|x|))^{-1}$, we have (recall that $h^{-1}$ is decreasing)
$$
h^{-1}(2r^{-1}) < |x| < h^{-1}(r^{-1}).
$$
Hence, we obtain
$$
I_{1} \leq C  \int_{(h(|x|))^{-1}}^{2(h(|x|))^{-1}}  h^{-1}(2r^{-1})^{-d-m} r t^{-\alpha-\beta}dr. 
$$
By \eqref{eqn 05.27.16:40}, we see that $h^{-1}(r^{-1}) \leq C h^{-1}(2r^{-1})$ for $r^{-1}<t_{1}^{-1}$. This and the above estimate for $I_1$ give 
$$
I_{1} \leq C  \int_{(h(|x|))^{-1}}^{2(h(|x|))^{-1}}  h^{-1}(r^{-1})^{-d-m} r t^{-\alpha-\beta}dr \leq \int_{(h(|x|))^{-1}}^{2t^{\alpha}}  h^{-1}(r^{-1})^{-d-m} r t^{-\alpha-\beta}dr.
$$
Note that for $r\geq (h(|x|))^{-1}$, we have $|x|\leq h^{-1}(r^{-1})$. This gives that
\begin{equation*}
\begin{aligned}
I_{2} &\leq C \int_{(h(|x|))^{-1}}^{t^{\alpha}}(h^{-1}(r^{-1}))^{-d-m} rt^{-\alpha-\beta}dr\\
& \leq C \int_{(h(|x|))^{-1}}^{2t^{\alpha}}(h^{-1}(r^{-1}))^{-d-m} rt^{-\alpha-\beta}dr.
\end{aligned}
\end{equation*}
For $I_{3}$, since $e^{-s} \leq c(\alpha) s^{1-\alpha}$ for $s\geq1$ and $s\mapsto (h^{-1}(s^{-1}))^{-d-m}$ is decreasing, we have (recall $|x|\leq h^{-1}(r^{-1})$)
\begin{equation*}
\begin{aligned}
I_{3} &\leq C \int_{t^{\alpha}}^{2t^{\alpha}} (h^{-1}(r^{-1}))^{-d-m} r t^{-\alpha-\beta} dr 
+C \int_{2t^{\alpha}}^{\infty} (h^{-1}(r^{-1}))^{-d-m} t^{-\beta} e^{-c(rt^{-\alpha})^{\frac{1}{1-\alpha}}}dr
\\
&\leq C \int_{(h(|x|))^{-1}}^{2t^{\alpha}} (h^{-1}(r^{-1}))^{-d-m} r t^{-\alpha-\beta} dr
 + C \left(h^{-1}\left(\frac{t^{-\alpha}}{2}\right)\right)^{-d-m}t^{\alpha-\beta}.
\end{aligned}
\end{equation*}
Using that $s\mapsto (h^{-1}(s^{-1}))^{-d-m}$ is decreasing again,
\begin{equation*}
\begin{aligned}
\left(h^{-1}\left(\frac{1}{2}t^{-\alpha}\right)\right)^{-d-m}t^{\alpha-\beta} &= C \int_{t^{\alpha}}^{2t^{\alpha}} (h^{-1}(\frac{t^{-\alpha}}{2}))^{-d-m} r t^{-\alpha-\beta} dr
\\
& \leq C \int_{t^{\alpha}}^{2t^{\alpha}} (h^{-1}(r^{-1}))^{-d-m} r t^{-\alpha-\beta} dr 
\\
& \leq C \int_{(h(|x|))^{-1}}^{2t^{\alpha}} (h^{-1}(r^{-1}))^{-d-m} r t^{-\alpha-\beta} dr.
\end{aligned}
\end{equation*} 
Thus, we have
$$
I_{1}+I_{2}+I_{3} \leq C \int_{(h(|x|))^{-1}}^{2t^{\alpha}} (h^{-1}(r^{-1}))^{-d-m} r t^{-\alpha-\beta} dr,
$$
and this proves the claim for $\beta\in\bN$ and $h(|x|)\leq t_{1}^{-1}$.

For $\beta\in\bN$ and $h(|x|)\geq t_{1}^{-1}$, if we use the relation $|x|\leq h^{-1}(r^{-1})$ for $r\geq (h(|x|))^{-1}$, we have
\begin{equation*}
\begin{aligned}
 \int_{t_1}^{\infty} |D^{m}_{x}p_{d} (r,x)| |\varphi_{\alpha,\beta} (t,r)| dr 
  \leq C\left(I_{4}+I_{3}\right),
\end{aligned}
\end{equation*}
where
\begin{gather*}
I_{4} = \int_{t_{1}}^{t^{\alpha}}  (h^{-1}(r^{-1}))^{-d-m} r t^{-\alpha-\beta} dr,
\end{gather*}
and $I_{3}$ comes from \eqref{eqn 05.28.12:51} which is bounded by
$$
C \int_{(h(|x|))^{-1}}^{2t^{\alpha}}  (h^{-1}(r^{-1}))^{-d-m} r t^{-\alpha-\beta} dr.
$$ 
Since $(h(|x|))^{-1} \leq t_{1}$, it is easy to see that
$$
I_{4} \leq C \int_{(h(|x|))^{-1}}^{2t^{\alpha}}  (h^{-1}(r^{-1}))^{-d-m} r t^{-\alpha-\beta} dr.
$$
Thus, we prove the lemma for the case $\beta\in\bN$. 

Finally, for the case $\beta\notin \bN$, by following all of the above computations with
$$
|\varphi_{\alpha,\beta}(t,r)| \leq C t^{-\beta} \quad\mbox{for}\;\; r\leq t^{\alpha},
$$
we have
\begin{align*}
\left|  \int_{t_1}^{\infty} |D^{m}_{x}p_{d} (r,x)| |\varphi_{\alpha,\beta} (t,r)| dr \right|
 \leq C \int_{(h(|x|))^{-1}}^{2t^{\alpha}} (h^{-1}(r^{-1}))^{-d-m}  t^{-\beta} dr.
\end{align*}
Thus, we obtain the first assertion. 

Now suppose that $\ell$ satisfies Assumption \ref{ell_con} (ii)--(1). Then using \eqref{hke_largetime2} we have
\begin{align*}
\int_{0}^{\infty} |D^{m}_{x}p_{d}(r,x)| |\varphi_{\alpha,\beta}(t,r)|dr \leq C \left( I'_{1} + I_{2}+ I_{3} \right),
\end{align*}
where
$$
I'_{1}=\int_{0}^{(h(|x|))^{-1}} |D^{m}_{x}p_{d}(r,x)| |\varphi_{\alpha,\beta}(t,r)| dr .
$$
Since $I'_{1}$ can be handled like $I_{1}$, the lemma is proved.
\end{proof}

\begin{corollary} \label{qintegrable}
Let $\alpha\in (0,1)$ and $\beta\in \bR$.

(i) There exists a constant $C=C(\alpha,\beta,d,\kappa_{1},\kappa_{2},\text{\boldmath{$\delta$}},\ell,C_{0},C_{1},C_{2})$ such that
\begin{equation*}
\int_{\bR^d} |q^{\alpha,\beta}_{d}(t,x)|dx \leq C t^{\alpha-\beta}, \quad t>0.
\end{equation*}

(ii) For any $0<\varepsilon<T<\infty$,
\begin{equation*}
\int_{\bR^d} \sup_{[\varepsilon,T]}  |q^{\alpha,\beta}_{d}(t,x)| dx < \infty.
\end{equation*}

\end{corollary}

\begin{proof}
(i) Due to the similarity, we only consider the case $\beta\in\bN$. 
 By Lemma \ref{lem 05.27.18:14},
\begin{align*}
\int_{\bR^d} |q^{\alpha,\beta}_d(t,x)|dx =&\int_{|x|\geq h^{-1}(t^{-\alpha})} |q^{\alpha,\beta}_d(t,x)|dx
+\int_{|x|< h^{-1}(t^{-\alpha})} |q^{\alpha,\beta}_d(t,x)|dx
\\
\leq&\; C \int_{|x|\geq h^{-1}(t^{-\alpha})} t^{2\alpha-\beta}\frac{K(|x|)}{|x|^d} dx\\
&+ C\int_{|x|< h^{-1}(t^{-\alpha})} \int_{0}^{t_1} p_{d} (r,x)\varphi_{\alpha,\beta} (t,r) dr dx \\
&+ C \int_{|x|< h^{-1}(t^{-\alpha})} \int_{(h(|x|))^{-1}}^{2t^{\alpha}}(h^{-1}(r^{-1}))^{-d} rt^{-\alpha-\beta} dr dx .
\end{align*}
Using that $\int_{\R^d}p_d(r,x)dx=1$ for all $r>0$ and $\int_{0}^{\infty} |\varphi_{\alpha,\beta}(t,r)|dr \leq C t^{\alpha-\beta}$, we have
$$
\int_{0}^{\infty}\int_{\bR^{d}} p_{d}(r,x) dx |\varphi_{\alpha,\beta}(t,r)| dr \leq C t^{\alpha-\beta}.
$$
Hence, by Fubini's theorem and \eqref{eqn 05.11.17:32}, we obtain
\begin{align*}
\int_{\bR^d} |q^{\alpha,\beta}_d(t,x)|dx
\leq&\; C \int_{r\geq h^{-1}(t^{-\alpha})} t^{2\alpha-\beta} \frac{K(r)}{r}dr + C t^{\alpha-\beta}
\\
&+ C \int_{0}^{2t^\alpha} \int_{(h(|x|))^{-1}\leq r} (h^{-1}(r^{-1}))^{-d} rt^{-\alpha-\beta} dx dr
\\
\leq&\; C t^{\alpha-\beta} + C \int_0^{2t^\alpha} rt^{-\alpha-\beta} dr
\leq C t^{\alpha-\beta}.
\end{align*}

(ii) Like (i), we only prove the case  $\beta\in\bN$. Let $0<\varepsilon<T<\infty$. Since $t^{2\alpha-\beta}\leq C(\varepsilon,T,\alpha,\beta)$ for $t\in[\varepsilon,T]$, by Lemma \ref{lem 05.27.18:14},
\begin{align*}
|q^{\alpha,\beta}_d(t,x)|\leq C(\alpha,\beta,d,\kappa_{1},\kappa_{2},\text{\boldmath{$\delta$}},\ell,C_{0},\varepsilon,T) \frac{K(|x|)}{|x|^d}, \quad t\in[\varepsilon,T].
\end{align*}
Also, if  $\varepsilon^\alpha h(|x|)\geq1$, and $t\in[\varepsilon,T]$, then   using Lemma \ref{lem 05.27.18:14} again, we obtain
\begin{align*}
|q^{\alpha,\beta}_d(t,x)|\leq C \left( \int_{0}^{t_1} p_{d} (r,x)\varphi_{\alpha,\beta} (t,r) dr+ \int_{(h(|x|))^{-1}}^{2T^\alpha} (h^{-1}(r^{-1}))^{-d} r dr \right).
\end{align*}
By Fubini's theorem and $\int_{\R^d}p_d(r,x)dx=1$ for all $r>0$, we have
\begin{equation*}
\begin{aligned}
&\int_{\bR^{d}}\sup_{t\in[\varepsilon,T]} \int_{0}^{t_1} p_{d} (r,x)\varphi_{\alpha,\beta} (t,r) dr dx\\
&\leq C \int_{0}^{\infty}\int_{\bR^{d}}  p_{d} (r,x) dx \left( \sup_{t\in[\varepsilon,T]} |\varphi_{\alpha,\beta}(t,r)| \right) dr
\\
&\leq C \int_{0}^{\infty}  \sup_{t\in[\varepsilon,T]} |\varphi_{\alpha,\beta}(t,r)|  dr
\\
&\leq C \int_{0}^{T^{\alpha}} r \varepsilon^{-\alpha-\beta}dr + C \int_{\varepsilon^{\alpha}}^{\infty} \varepsilon^{-\beta} e^{-c(rT^{-\alpha})^{\frac{1}{1-\alpha}}} dr \leq C.
\end{aligned}
\end{equation*}
Hence, like the proof of (i),
\begin{align*}
&\int_{\bR^d} \sup_{[\varepsilon,T]}|q^{\alpha,\beta}_{d}(t,x)|dx\\
& = \int_{|x|\geq h^{-1}(\varepsilon^{-\alpha})} \sup_{[\varepsilon,T]}|q^{\alpha,\beta}_{d}(t,x)|dx
+\int_{|x|< h^{-1}(\varepsilon^{-\alpha})} \sup_{[\varepsilon,T]}|q^{\alpha,\beta}_{d}(t,x)|dx
\\
&\leq C \int_{|x|\geq h^{-1}(\varepsilon^{-\alpha})} \frac{K(|x|)}{|x|^d} dx + C \int_{\bR^{d}} \sup_{[\varepsilon,T]} \int_{0}^{t_1} p_{d} (r,x)\varphi_{\alpha,\beta} (t,r) drdx 
\\
&\quad+ C \int_{|x|< h^{-1}(\varepsilon^{-\alpha})} \int_{(h(|x|))^{-1}}^{2T^{\alpha}}(h^{-1}(r^{-1}))^{-d} r dr dx
\\
&\leq C + C \int_{0}^{2T^\alpha} \int_{(h(|x|))^{-1}\leq r} (h^{-1}(r^{-1}))^{-d} r dx dr
\leq C + C \int_0^{2T^\alpha} r dr <\infty.
\end{align*}
The corollary is proved.
\end{proof}

Recall that $C_{p}^{\infty}(\bR^d)$ is the set of functions $u_0=u_0(x)$ such that $D^m_xu_0\in L_p$ for any $m\in \bN_{0}$.  
\begin{lemma} \label{zero converge}
Let $u_0\in C_{p}^{\infty}(\bR^d)$, and define $u$ as
$$
u(t,x):=\int_{\bR^d} q_{d}(t,x-y)u_0(y)dy.
$$

(i) As as $t \downarrow 0$, $u(t,\cdot)$ converges to $u_0(\cdot)$ uniformly on $\bR^d$ and also  in $H^n_p$ for any $n\in \bN_0$.

(ii) $u\in C^{\alpha,\infty}_{p}([0,T]\times\bR^d)$ and $u$ satisfies \eqref{eqn 06.15.16:05}. 

\end{lemma}
\begin{proof}
We can prove the lemma by following proof of \cite[Lemma 5.1]{kim2020nonlocal}.
\end{proof}

\medskip

\mysection{Estimation of solution: Calder\'on-Zygmund approach}\label{sec 5}

In this section we prove some a priori estimates for solutions to the  equation with zero initial condition
\begin{equation}\label{mainequation-1}
\partial^{\alpha}_{t}u = \mathcal{L}u + f,\quad t>0\,; \quad u(0,\cdot)=0.
\end{equation}
We first provide the representation formula. 

\begin{lemma} \label{u=qfsolution}

(i) Let $u\in C_c^\infty(\bR_+^{d+1})$ and denote $f:=\partial^{\alpha}_{t}u-\mathcal{L}u$. Then
\begin{equation} \label{u=qf}
\begin{gathered}
u(t,x)=\int_{0}^{t} \int_{\bR^{d}} q^{\alpha,1}_{d}(t-s,x-y) f(s,y) dy ds.
\end{gathered}
\end{equation}

(ii) Let $f\in C_c^\infty(\bR_+^{d+1})$ and define $u$ as in \eqref{u=qf}. Then $u$ satisfies equation \eqref{mainequation-1} for each $(t,x)$.
\end{lemma}

\begin{proof}
The lemma is an extension of \cite[Lemma 4.1]{kim2020nonlocal} which treats the case $\psi(|\xi|)=\phi(|\xi|^{2})$, where $\phi$ is a Bernstein function. The proof of \cite[Lemma 4.1]{kim2020nonlocal} only uses \eqref{fourier of q} with $\phi(|\xi|^{2})$ in place of $\psi(|\xi|)$ and Corollary \ref{qintegrable}. Since no property of $\phi$ is used, we can prove all the claims by repeating the same argument. For more detail, see also \cite[Lemma 3.5]{kim17timefractionalpde}.
\end{proof}

For $f\in \Ccinf(\bR^{d+1})$, we define 
\begin{align}
L_{0}f(t,x) := \int_{-\infty}^{t}\int_{\bR^{d}} q^{\alpha,1}_{d}(t-s,x-y)f(s,y) dy ds,
\nn\\
Lf(t,x) := \int_{-\infty}^{t}\int_{\bR^{d}} q^{\alpha,1}_{d}(t-s,x-y)\mathcal{L}f(s,y) dy ds,\label{def_Lf}
\end{align}
where $\mathcal{L}f(s,y):= \mathcal{L}(f(s,\cdot))(y)$.

Note that  
$\mathcal{L}f$ is bounded for any $f\in C_c^\infty (\bR^{d+1})$. Thus, the operator $L$ is well defined on $C_c^\infty (\bR^{d+1})$.
For each fixed $s$ and $t$ such that  $s<t$, define
\begin{align*}
\Lambda^{0}_{t,s}f(x) &:=\int_{\bR^d} q^{\alpha,1}_{d}(t-s,x-y) \mathcal{L} f(s,y) dy,\\
\Lambda_{t,s}f(x) &:=\int_{\bR^d} q^{\alpha,\alpha+1}_{d}(t-s,x-y)  f(s,y) dy.
\end{align*}
Note that, by Corollary \ref{qintegrable}, $\Lambda_{t,s}f$ and $\Lambda^{0}_{t,s}f$ are square integrable. Moreover, using \eqref{fourier of q} and the recurrence relation 
$E_{a,b}(z)=\frac{1}{\Gamma(b)} + zE_{a,a+b}(z)$ ($z\in \bC$)
we have
\begin{equation*}
\begin{aligned}
\cF_{d}\{q^{\alpha,\alpha+1}_{d}(t-s,\cdot)\}(\xi) &= (t-s)^{-1}E_{\alpha,0}(-(t-s)^{\alpha}\psi(|\xi|))
\\
& = -(t-s)^{\alpha-1}\psi(|\xi|) E_{\alpha,\alpha}(-(t-s)^{\alpha}\psi(|\xi|)) 
\\
&= -\psi(|\xi|) \cF_{d}\{q^{\alpha,1}_{d}(t-s,\cdot)\}(\xi).
\end{aligned}
\end{equation*}
Hence,
\begin{equation*}\label{eqn 06.12.15:17}
\begin{aligned}
\cF_d \{\Lambda^{0}_{t,s}f\}(\xi) & = -\psi(|\xi|)  \cF_{d}\{q^{\alpha,1}_{d}\}(t-s,\xi) \hat{f}(s,\xi)
\\
&=\mathcal{F}_{d}\{q^{\alpha,\alpha+1}_{d} \}(t-s,\xi)\hat{f}(s,\xi)=\cF_d \{\Lambda_{t,s}f\}(\xi).
\end{aligned}
\end{equation*}
Thus, we have 
\begin{equation}\label{eqn 06.02.15:23}
\begin{aligned}
\mathcal{L}L_{0}f(t,x) &= Lf(t,x)
= \lim_{\varepsilon \downarrow 0} \int_{-\infty}^{t-\varepsilon}\left(\int_{\bR^{d}} q^{\alpha,\alpha+1}_{d}(t-s,x-y) f(s,y) dy \right) ds.
\end{aligned}
\end{equation}

We give a short description of our strategy to obtain the estimation of solution. By \eqref{u=qf}, and \eqref{eqn 06.02.15:23} if we let $u$ be a solution to equation \eqref{mainequation-1}, then, we have $\mathcal{L}u=Lf$. Hence to obtain the estimation of solution, it is enough to show the $L_{q}(L_{p})$-boundedness of the linear operator  $L$ (Theorem \ref{thm 05.14.18:28}).

We prove Theorem \ref{thm 05.14.18:28} step by step by using the following lemmas;
\begin{enumerate}[(i)]
\item
We prove $L_{2}$-boundedness of $L$ in Lemma \ref{22estimate};

\item
From Lemma \ref{ininestimate} to Lemma \ref{outinspaceestimate}, we control the mean oscillation of $L$ in terms of $\|f\|_{L_{\infty}(\bR^{d+1})}$;

\item
We prove Theorem \ref{thm 05.14.18:28} by using the result of (ii) and Theorem \ref{feffermanstein}.
\end{enumerate}

\begin{remark}\label{rmk 10.10.17:10}
For parabolic (i.e. $\alpha=1$) equations, we see from \eqref{eqn 06.02.15:23} that the above argument requires estimation of $\mathcal{L}p_{d}$. To obtain estimation on $\mathcal{L}p_{d}$ we need more differentiability on $j_{d}$ and corresponding conditions in Assumption \ref{ass bernstein} (i), (ii). Thus, to cover more general $j_{d}$, we only consider the case $\alpha\in(0,1)$ in this article.
\end{remark}

The following two lemmas will be used for step (ii). 
For notational convenience, we define a strictly increasing function $\kappa:(0,\infty)\to (0,\infty)$ as $\kappa(r) := (h(r))^{-1/\alpha}$.

Recall that we assume that Assumption \ref{ass bernstein}, Assumption \ref{ass bernstein-2} and Assumption \ref{ell_con}.
\begin{lemma}\label{l:qest_case1} 
Let $\alpha\in(0,1)$. There exists $C=C(\alpha,d,\kappa_{1},\kappa_{2},\text{\boldmath{$\delta$}},\ell,C_{0},C_{1},C_{2})$ such that for any $b>0$
\begin{align}
&\int_{\kappa(b)}^\infty \int_{|y|\ge b } |D_xq_d^{\alpha,\alpha+1}(s,y)| dy ds\le Cb^{-1},\label{q1upper}\\
&\int_{\kappa(4b)}^{\infty} \int_{|y|\leq 4b} |q_d^{\alpha,\alpha+1}(s,y)| dy ds\le C.\label{q0upper}
\end{align}
\end{lemma}

\begin{proof}
Let $t_1>0$ be the constant in Lemma \ref{l:hke_largetime}. We will consider two cases separately.

\noindent{\bf (Case 1)} $\ell$ satisfies Assumption \ref{ell_con} (i).

We first show \eqref{q1upper}. By \eqref{eqn 07.02.16:48}, 
\begin{equation}\label{eqn 09.06.19:51}
\begin{aligned}
&\int_{\kappa(b)}^\infty \int_{|y| \geq b} |D_xq_d^{\alpha,\alpha+1}(s,y)| dy ds\\
&\le\int_{\kappa(b)}^\infty \int_{b\leq |y| \leq h^{-1}(s^{-\alpha})} \int_{0}^{\infty} |D_{x}p_{d} (r,y)| |\varphi_{\alpha,\alpha+1} (s,r)| dr dy ds\\
&\quad+C\int_{\kappa(b)}^\infty \int_{|y| \geq h^{-1}(s^{-\alpha})} s^{\alpha-1}\frac{K(|y|)}{|y|^{d+1}} dy ds\\
&=\int_{\kappa(b)}^\infty \int_{b\leq |y| \leq h^{-1}(s^{-\alpha})} \int_{0}^{t_1} |D_{x}p_{d} (r,y)| |\varphi_{\alpha,\alpha+1} (s,r)| dr dy ds\\
&\quad+\int_{\kappa(b)}^\infty \int_{b\leq |y| \leq h^{-1}(s^{-\alpha})} \int_{t_1}^{\infty} |D_{x}p_{d} (r,y)| |\varphi_{\alpha,\alpha+1} (s,r)| dr dy ds\\
&\quad+C\int_{\kappa(b)}^\infty \int_{|y| \geq h^{-1}(s^{-\alpha})} s^{\alpha-1}\frac{K(|y|)}{|y|^{d+1}} dy ds\\
&=:I_1+I_2+ C I_3.
\end{aligned}
\end{equation}
By \eqref{q1upper_case1_pf}, Lemma \ref{l:qest_largetime} and Lemma \ref{l:qest_int}, we obtain $I_1+I_2\le Cb^{-1}$.  
By \eqref{eqn 05.11.17:32} and the change of variable ($s^{\alpha}\to s$),
\begin{equation}\label{eqn 08.30.11:05}
\begin{aligned}
I_3 & =\int_{(h(b))^{-1/\alpha}}^{\infty} \int_{h^{-1}(s^{-\alpha})}^{\infty} s^{\alpha-1}\frac{K(\rho)}{\rho^{2}} d\rho ds
\\
& \leq C\int_{(h(b))^{-1/\alpha}}^{\infty} h^{-1}(s^{-\alpha})^{-1}s^{\alpha-1}\int_{h^{-1}(s^{-\alpha})}^{\infty} \frac{K(\rho)}{\rho} d\rho ds
\\
& \leq C\int_{(h(b))^{-1/\alpha}}^{\infty} h^{-1}(s^{-\alpha})^{-1}s^{-1}ds
\\
& = C\int_{(h(b))^{-1}}^{\infty} h^{-1}(s^{-1})^{-1}s^{-1}ds.
\end{aligned}
\end{equation}
Thus, by Lemma \ref{l:qest_offd} with $f(r)=h(r^{-1})$ and $k=1$, we obtain $I_3\le Cb^{-1}$.

Now, we show \eqref{q0upper}. By \eqref{q0upper_case1_pf} and Lemma \ref{l:qest_largetime},
\begin{equation}\label{eqn 09.06.20:00}
\begin{aligned}
&\int_{\kappa(4b)}^{\infty} \int_{|y|\leq 4b} q_d^{\alpha,\alpha+1}(s,y) dy ds\\
&\le \int_{\kappa(4b)}^\infty \int_{|y| \leq 4b} \int_{0}^{t_1} p_{d} (r,y) |\varphi_{\alpha,\alpha+1} (s,r)| dr dy ds\\
&\quad+ \int_{\kappa(4b)}^\infty \int_{|y| \leq 4b} \int_{t_1}^{\infty} p_{d} (r,y) |\varphi_{\alpha,\alpha+1} (s,r)| dr dy ds\\
&\le C+ C\int_{\kappa(4b)}^\infty \int_{|y| \leq 4b} \int_{(h(|y|))^{-1}}^{2s^{\alpha}} (h^{-1}(r^{-1}))^{-d}  s^{-\alpha-1} dr dy ds
\\
& = C +CI_{4}.
\end{aligned}
\end{equation}
Note that if $\kappa(4b)<s$, then by Fubini's theorem
\begin{equation*}
\begin{aligned}
I_{4}&=\int_{\kappa(4b)}^{\infty} \int_{0}^{4b} \int_{(h(\rho))^{-1}}^{2s^{\alpha}} (h^{-1}(r^{-1}))^{-d} s^{-\alpha-1} \rho^{d-1} dr  d\rho ds
\\
& = \int_{(h(4b))^{-1/\alpha}}^{\infty} \int_{0}^{4b} \int_{(h(\rho))^{-1}}^{(h(4b))^{-1}} (h^{-1}(r^{-1}))^{-d} s^{-\alpha-1} \rho^{d-1} dr  d\rho ds
\\
&  \quad + \int_{(h(4b))^{-1/\alpha}}^{\infty} \int_{0}^{4b} \int_{(h(4b))^{-1}}^{2s^{\alpha}} (h^{-1}(r^{-1}))^{-d} s^{-\alpha-1} \rho^{d-1} dr  d\rho ds
\\
&\leq C \int_{(h(4b))^{-1/\alpha}}^{\infty} \int_{0}^{(h(4b))^{-1}} \int_{0}^{h^{-1}(r^{-1})} (h^{-1}(r^{-1}))^{-d} s^{-\alpha-1} \rho^{d-1}d\rho  dr  ds
\\
& \quad  + C b^{d} \int_{(h(4b))^{-1/\alpha}}^{\infty}\int_{(h(4b))^{-1}}^{2s^{\alpha}} (h^{-1}(r^{-1}))^{-d} s^{-\alpha-1}  dr  ds.
\end{aligned}
\end{equation*}
We can easily check
\begin{equation}\label{eqn 08.30.12:00}
\begin{aligned}
& \int_{(h(4b))^{-1/\alpha}}^{\infty} \int_{0}^{(h(4b))^{-1}} \int_{0}^{h^{-1}(r^{-1})} (h^{-1}(r^{-1}))^{-d} s^{-\alpha-1} \rho^{d-1}  d\rho dr ds 
\\
& \leq \int_{(h(4b))^{-1/\alpha}}^{\infty} \int_{0}^{(h(4b))^{-1}}  s^{-\alpha-1} dr ds \leq C. 
\end{aligned}
\end{equation}
By Fubini's theorem and Lemma \ref{l:qest_offd} with $f(r)=h(r^{-1})$ and $k=d$, we have
\begin{equation}\label{eqn 08.30.12:00-2}
\begin{aligned}
&b^{d} \int_{(h(4b))^{-1/\alpha}}^{\infty}\int_{(h(4b))^{-1}}^{2s^{\alpha}} (h^{-1}(r^{-1}))^{-d} s^{-\alpha-1}  dr  ds
\\
&=b^{d} \int_{(h(4b))^{-1}}^{\infty}\int_{(r/2)^{1/\alpha}} s^{-\alpha-1} (h^{-1}(r^{-1}))^{-d} dsdr
\\
& \leq C b^{d} \int_{(h(4b))^{-1}}^{\infty}r^{-1} (h^{-1}(r^{-1}))^{-d} dr=C.
\end{aligned}
\end{equation}

\noindent{\bf (Case 2)} $\ell$ satisfies the Assumption \ref{ell_con} (ii).

Suppose that $\ell$ satisfies Assumption \ref{ell_con} (ii)--(2).
We first show \eqref{q1upper}.  Using \eqref{eqn 09.03.19:06} and \eqref{eqn 07.02.16:48}
\begin{align*}
&\int_{\kappa(b)}^\infty \int_{|y| \geq b} |D_xq_d^{\alpha,\alpha+1}(s,y)| dy ds\\
&\le\int_{\kappa(b)}^\infty \int_{b\leq |y| \leq h^{-1}(s^{-\alpha})} \int_{0}^{\infty} |D_{x}p_{d} (r,y)| |\varphi_{\alpha,\alpha+1} (s,r)| dr dy ds\\
&\quad+C\int_{\kappa(b)}^\infty \int_{|y| \geq h^{-1}(s^{-\alpha})} s^{\alpha-1}\frac{K(|y|)}{|y|^{d+1}} dy ds\\
&=\int_{\kappa(b)}^\infty \int_{b\leq |y| \leq h^{-1}(s^{-\alpha})} \int_{0}^{t_1} |D_{x}p_{d} (r,y)| |\varphi_{\alpha,\alpha+1} (s,r)| dr dy ds\\
&\quad+\int_{\kappa(b)}^\infty \int_{b\leq |y| \leq h^{-1}(s^{-\alpha})} \int_{t_1}^{\infty} |D_{x}p_{d} (r,y)| |\varphi_{\alpha,\alpha+1} (s,r)| dr dy ds\\
&\quad+C\int_{\kappa(b)}^\infty \int_{|y| \geq h^{-1}(s^{-\alpha})} s^{\alpha-1}\frac{K(|y|)}{|y|^{d+1}} dy ds\\
&=: I_5+I_6+ C I_3 .
\end{align*}
Using \eqref{eqn 08.30.11:05} and \eqref{q1upper_case2_pf} we have $I_{3}+I_{5} \leq Cb^{-1}$. Also, by Lemma \ref{l:qest_largetime} and Lemma \ref{l:qest_int} with $f(r)=(h(r))^{-1}$, we have $I_{6}\leq Cb^{-1}$, and thus we have desired result.

To show \eqref{q0upper}, we use \eqref{q0upper_case2_pf} and  Lemma \ref{l:qest_largetime} for  \eqref{eqn 09.06.20:00}. Then, by \eqref{eqn 08.30.12:00} and \eqref{eqn 08.30.12:00-2}, we have
\begin{align*}
\int_{\kappa(4b)}^{\infty} \int_{|y|\leq 4b} q_d^{\alpha,\alpha+1}(s,y) dy ds \leq C + CI_{4} \leq C.
\end{align*}

 Now suppose that $\ell$ satisfies Assumption \ref{ell_con} (ii)--(1). Then by Lemma \ref{l:qest_largetime}
 for all $(t,x)\in(0,\infty)\times\bR^{d}$ such that $t^{\alpha}h(|x|)\geq 1$ we have,
\begin{align*}
& |D^{m}_{x}q^{\alpha,\alpha+1}_{d}(t,x)|  \leq C  \int_{(h(|x|))^{-1}}^{2t^{\alpha}} (h^{-1}(r^{-1}))^{-d-m}  t^{-\alpha-1} dr.
\end{align*}
This implies that when we show \eqref{q1upper}, from \eqref{eqn 09.06.19:51}, we only need to handle $I_{3}$ and
$$
\int_{\kappa(b)}^{\infty}\int_{b\leq |y|\leq h^{-1}(s^{-\alpha})} \int_{(h(|y|))^{-1}}^{2s^{\alpha}} (h^{-1}(r^{-1}))^{-d-1}  s^{-\alpha-1} dr.
$$
Hence, using Lemma \ref{l:qest_int} with $f(r)=h(r^{-1})$ and \eqref{eqn 08.30.11:05}, we have \eqref{q1upper}.
For \eqref{q0upper}, we can check that from \eqref{eqn 09.06.20:00}, we only need to handle $I_{4}$. Since $I_{4}\leq C$, we prove all claims. The lemma is proved.
\end{proof}

Recall that the operator $L$ defined in \eqref{def_Lf} satisfies \eqref{eqn 06.02.15:23}.
 
\begin{lemma} \label{22estimate}
For any $f\in\Ccinf(\bR^{d+1})$, we have
\begin{equation*}
\begin{gathered}
\|L f\|_{L_2(\bR^{d+1})}\leq C(\alpha,d) \|f\|_{L_2(\bR^{d+1})}.
\end{gathered}
\end{equation*}
Consequently, the  operators $L$ can be  continuously extended to $L_2(\bR^{d+1})$. 
\end{lemma}

\begin{proof}
Like Lemma \ref{u=qfsolution}, we can prove the lemma by following the proof of \cite[Lemma 4.2]{kim2020nonlocal}. Again note that the only difference between our lemma and \cite[Lemma 4.2]{kim2020nonlocal} is that we just replace $\psi(|\xi|)$ in place of $\phi(|\xi|^{2})$ since we do not need any property of $\phi$.
\end{proof}

Recall that $\kappa(b) = (h(b))^{-1/\alpha}.$
For $(t,x)\in\R^{d+1}$ and $b>0$,  denote
\begin{equation*}
 Q_b(t,x)=(t-\kappa(b),\,t+\kappa(b))\times {B}_b(x),
\end{equation*}
and
\begin{equation*}
Q_b=Q_b(0,0), \quad B_b=B_b(0).
\end{equation*}
 For  measurable subsets $Q\subset \R^{d+1}$ with  finite measure and  locally integrable functions $f$, define
\begin{equation*}
f_Q=\aint_{Q}f(s,y)dyds=\frac{1}{|Q|}\int_{Q}f(s,y)dyds,
\end{equation*}
where $|Q|$ is the Lebesgue measure of $Q$.

\begin{lemma} \label{ininestimate}
Let $b>0$, and $f\in C^{\infty}_c(\bR^{d+1})$  have a support in   $ (-3\kappa(b), 3\kappa(b))\times B_{3b}$. Then,
\begin{equation*}
\begin{gathered}
\aint_{Q_b}|L f (t,x)|dxdt \leq C(\alpha,d) \|f\|_{L_\infty(\bR^{d+1})}.
\end{gathered}
\end{equation*}
\end{lemma}

\begin{proof}
By H\"older's inequality and Lemma \ref{22estimate},
\begin{align*}
\aint_{Q_b}|L f (t,x)|dxdt \leq &  \;\|L f\|_{L_2(\bR^{d+1})}|Q_b|^{-1/2}
\leq  \;C \|f\|_{L_2(\bR^{d+1})}|Q_b|^{-1/2}
\\
= & \;\left( \int_{-3\kappa(b)}^{3\kappa(b)} \int_{B_{3b}} |f(t,x)|^{2} dydt \right)^{1/2}|Q_b|^{-1/2}
\leq  C \|f\|_{L_\infty(\bR^{d+1})}.
\end{align*}
The lemma is proved.
\end{proof}

\begin{lemma} \label{inwholeestimate} 
Let $b>0$, and $f\in C_c^\infty(\bR^{d+1})$  have a support in $(-3\kappa(b), \infty)\times \bR^d$. Then,
\begin{equation*}
\begin{gathered}
\aint_{Q_b}|L f (t,x)|dxdt \leq C(\alpha,d,\kappa_{1},\kappa_{2},\text{\boldmath{$\delta$}},\ell,C_{0},C_{1},C_{2}) \|f\|_{L_\infty(\bR^{d+1})}.
\end{gathered}
\end{equation*}
\end{lemma}

\begin{proof}
Take $\eta=\eta(t)\in C^{\infty}(\bR)$ such that $0\leq \eta\leq 1$, $\eta(t)=1$ for $t\leq 2\kappa(b)$, and $\eta(t)=0$ for $t\geq 5\kappa(b)/2$.
Since $L f=L (f\eta)$ on $Q_b$ and $|f\eta|\leq |f|$, it is enough to assume $f(t,x)=0$ if $|t|\geq 3\kappa(b)$ to prove the lemma.

Choose a function $\zeta=\zeta(x) \in C_c^\infty (\bR^d)$ such that $\zeta=1$ in $B_{7b/3}$, $\zeta=0$ outside of $B_{8b/3}$ and $0\leq\zeta\leq1$. Set $f_1=\zeta f$ and $f_2=(1-\zeta)f$. Then due to the linearity, $L f = L f_1 + L f_2$. Also, since $Lf_{1}$ can be estimated by Lemma \ref{ininestimate}, to prove the lemma, we may further assume that $f(t,y)=0$ if $y\in B_{2b}$. Hence, for any $x\in B_b$,
\begin{align*}
\int_{\bR^d} \left|q^{\alpha,\alpha+1}_{d} (t-s,x-y) f(s,y)\right| dy =& \int_{|y-x|\geq 2b} \left|q^{\alpha,\alpha+1}_{d} (t-s,y) f(s,x-y)\right| dy 
\\
\leq & \int_{|y|\geq b} \left|q^{\alpha,\alpha+1}_{d} (t-s,y) f(s,x-y) \right| dy.
\end{align*}
Using Lemma \ref{lem 05.27.18:14} and \eqref{eqn 05.11.17:32},
\begin{align*}
& \int_{|y|\geq b} \left|q^{\alpha,\alpha+1}_{d} (t-s,y) f(s,x-y) \right| dy 
\\
&\leq \|f\|_{L_\infty(\bR^{d+1})} 1_{|s|\leq 3\kappa(b)} \int_{|y|\geq b} |q^{\alpha,\alpha+1}_{d} (t-s,y)| dy
\\
&\leq C\|f\|_{L_\infty(\bR^{d+1})}1_{|s|\leq 3\kappa(b)} \int_b^\infty (t-s)^{\alpha-1} \frac{K(\rho)}{\rho^{d}}\rho^{d-1}d\rho
\\
&\leq C\|f\|_{L_\infty(\bR^{d+1})} 1_{|s|\leq 3\kappa(b)}  (t-s)^{\alpha-1} h(b).
\end{align*}
Note that if $|t| \leq \kappa(b)$ and $|s|\leq 3\kappa(b)$ then $|t-s|\leq 4\kappa(b)$. Thus we have
\begin{align*}
|L f(t,x)| &\leq C\|f\|_{L_\infty(\bR^{d+1})} h(b) \int_{|t-s|\leq 4\kappa(b)} |t-s|^{\alpha-1} ds
 \leq C \|f\|_{L_{\infty}(\R^{d+1})}.
\end{align*}
This implies the desired estimate. The lemma is proved.
\end{proof}

\begin{lemma}
\label{outwholetimeestimate}
Let $b>0$, and $f\in C_c^\infty(\bR^{d+1})$  have a support in  $(-\infty, -2\kappa(b))\times \bR^d$. Then there is $C=C(\alpha,d,\kappa_{1},\kappa_{2},\text{\boldmath{$\delta$}},\ell,C_{0},C_{1},C_{2})$ such that for  any $(t_1,x), (t_2,x)\in Q_b$, 
\begin{gather*}
\aint_{Q_b}\aint_{Q_b} |L f (t_1,x)-L f(t_2,x)| dx dt_1 d\tilde{x} dt_2\leq C \|f\|_{L_\infty(\bR^{d+1})}.
\end{gather*}

 \end{lemma}

\begin{proof}
We will show that 
$$
|L f (t_1,x)-L f(t_2,x)|   \leq C \|f\|_{L_\infty(\bR^{d+1})},
$$
and this certainly proves the lemma. 

Without loss of generality, we assume that $t_{1}>t_{2}$. Recall $f(s,y)=0$ if $s\geq -2\kappa(b)$. Thus,  if $t>-\kappa(b)$,  by  applying the fundamental theorem of calculus and \eqref{eqn 07.28 20:25}, we have
\begin{align*}
|L f(t_1,x)-L f(t_2,x)|
&= \left|\int_{-\infty}^{-2\kappa(b)}  \int_{\bR^d}  \int_{t_2}^{t_1} q^{\alpha,\alpha+2}_{d} (t-s,x-y) f(s,y)  dt dy ds \right|.
\end{align*}
Using Corollary \ref{qintegrable} (i),
\begin{align*}
&\int_{\bR^d} \int_{t_2}^{t_1} |q^{\alpha,\alpha+2}_{d} (t-s,x-y) f(s,y)| dt dy
\\
&\leq C \|f\|_{L_{\infty}} \int_{t_{2}}^{t_{1}} \int_{\bR^{d}} |q^{\alpha,\alpha+2}_{d}(t-s,y)| dy dt
\leq C \|f\|_{L_\infty(\bR^{d+1})} \int_{t_2}^{t_1} (t-s)^{-2}dt.
\end{align*}
Thus,  if  $-\kappa(b)\leq t_2<t_1 \leq \kappa(b)$, 
\begin{align*}
&\Big|\int_{-\infty}^{-2\kappa(b)}  \int_{\bR^d}  \int_{t_2}^{t_1} q^{\alpha,\alpha+2}_{d} (t-s,x-y) f(s,y)  dt dy ds\Big|
\\
&\leq C \|f\|_{L_\infty(\bR^{d+1})} \left(\int_{t_2}^{t_1}\int_{-\infty}^{-2\kappa(b)} (t-s)^{-2}ds dt \right)
\\
&\leq C \|f\|_{L_\infty(\bR^{d+1})} \left(\int_{t_2}^{t_1} \frac{1}{\kappa(b)} dt \right)
\leq C \|f\|_{L_\infty(\bR^{d+1})}.
\end{align*}
This completes the proof.
\end{proof}

\begin{lemma}
\label{outoutspaceestimate}
Let $b>0$, and $f\in C_c^\infty(\bR^{d+1})$ have a support in  $(-\infty,-2\kappa(b))\times B_{2b}^c$. Then there is $C=C(\alpha,d,\kappa_{1},\kappa_{2},\text{\boldmath{$\delta$}},\ell,C_{0},C_{1},C_{2})$ such that for any $(t,x_{1}),(t,x_{2})\in Q_b$, 
$$
\aint_{Q_{b}}\aint_{Q_b} |L f (t,x_{1}) - L f(t,x_{2})| dx_{1} dt dx_{2} d\tilde{t} \leq C(\alpha,d,\kappa_{1},\kappa_{2},\text{\boldmath{$\delta$}},\ell,C_{0}) \|f\|_{L_\infty(\bR^{d+1})}.
$$
\end{lemma}

\begin{proof}
Like the previous lemma it is enough to show that
\begin{gather*}
|L f (t,x_{1}) - L f(t,x_{2})| \leq C \|f\|_{L_\infty(\bR^{d+1})}.
\end{gather*}

Recall $f(s,y)=0$ if $s\geq -2\kappa(b)$ or $|y|\leq 2b$. Thus,  if $t>-\kappa(b)$, by the fundamental theorem of calculus,
\begin{equation*}
\begin{aligned}
&|L f(t,x_{1}) - L f(t,x_{2})|
\\
&=\Big|\int_{-\infty}^{-2\kappa(b)}  \int_{|y|\geq 2b} \left(q^{\alpha,\alpha+1}_{d} (t-s,x_{1}-y) -q^{\alpha,\alpha+1}_{d}(t-s,x_{2}-y) \right) f(s,y) dy ds \Big|
\\
&= \Big|\int_{-\infty}^{-2\kappa(b)}  \int_{|y|\geq 2b} \int_{0}^{1}  \nabla_{x} q^{\alpha,\alpha+1}_{d}(t-s,\lambda(x_{1},x_{2},u)-y) \cdot (x_{1}-x_{2}) du f(s,y) dy ds \Big|,
\end{aligned}
\end{equation*}
where $\lambda(x_{1},x_{2},u)=ux_{1}+(1-u)x_{2}$. Since $x_{1},x_{2}\in B_{b}$ and $|y|\geq 2b$, we can check that $|\lambda(x_{1},x_{2},u)-y| \geq b$. Thus, by the change of variable $(\lambda(x_{1},x_{2},u)-y) \to y$,
\begin{equation*}
\begin{aligned}
|L f(t,x_{1}) - L f(t,x_{2})|
& \leq C  b \|f\|_{L_\infty(\bR^{d+1})}  \int_{-\infty}^{ -2\kappa(b)} \int_{|y|\geq b} |\nabla_{x}q^{\alpha,\alpha+1}_{d} (t-s,y)| dy ds
\\
& \leq C b  \|f\|_{L_\infty(\bR^{d+1})}  \int_{\kappa(b)}^{\infty} \int_{|y|\geq b} |\nabla_{x}q^{\alpha,\alpha+1}_{d} (s,y)| dy ds.
\end{aligned}
\end{equation*}
Thus, by \eqref{q1upper}, we obtain the desired result.
\end{proof}

\begin{lemma} \label{outinspaceestimate}
Let $b>0$, and $f\in C_c^\infty(\bR^{d+1})$  have a support in   $(-\infty,-2\kappa(b))\times B_{3b}$. Then for any $(t,x)\in Q_b$
\begin{gather*}
\aint_{Q_b}|L f (t,x)|dx dt \leq C(\alpha,d,\kappa_{1},\kappa_{2},\text{\boldmath{$\delta$}},\ell,C_{0},C_{1},C_{2})\|f\|_{L_\infty(\bR^{d+1})}.
\end{gather*} 
\end{lemma}

\begin{proof}
For $(t,x)\in Q_b$,
\begin{align*}
|L f(t,x)| \leq&\int_{-\infty}^{-2\kappa(b)} \int_{B_{3b}} |q^{\alpha,\alpha+1}_{d} (t-s,x-y)f(s,y)|dyds
\\
\leq& C \|f\|_{L_\infty} \int_{-\infty}^{-2\kappa(b)} \int_{B_{3b}} |q^{\alpha,\alpha+1}_{d} (t-s,x-y)|dyds
\\
\leq& C \|f\|_{L_\infty} \int_{\kappa(b)}^{\infty} \int_{B_{4b}} |q^{\alpha,\alpha+1}_{d} (s,y)|dyds =: C \|f\|_{L_\infty} \left( I +II \right),
\end{align*}
where
$$
I=\int^{\kappa(4b)}_{\kappa(b)} \int_{B_{4b}} |q^{\alpha,\alpha+1}_{d} (s,y)|dyds, \quad  \quad
II=\int_{\kappa(4b)}^{\infty} \int_{B_{4b}} |q^{\alpha,\alpha+1}_{d} (s,y)|dyds.
$$
By Corollary \ref{qintegrable} (i) and the scaling properties of $h$ or $\ell$ 
\begin{align*}
I &= \int^{\kappa(4b)}_{\kappa(b)} \int_{B_{4b}}  |q^{\alpha,\alpha+1}_{d} (s,y)|dyds
\leq C \int^{\kappa(4b)}_{\kappa(b)} s^{-1}ds
 \leq C \log\left(\frac{\kappa(4b)}{\kappa(b)} \right) 
\leq C \log\left(16\right).
\end{align*}
By \eqref{q0upper}, we also have $II\le C$.
Thus, $I$ and $II$ are bounded by a constant independent of $b$. Hence, we have
$$
|L f(t,x)|\leq C  \|f\|_{L_\infty(\bR^{d+1})},
$$ 
and the lemma is proved.
\end{proof}

\begin{corollary}
\label{outwholeestimate}
There is $C=C(\alpha,d,\kappa_{1},\kappa_{2},\text{\boldmath{$\delta$}},\ell,C_{0},C_{1},C_{2})$ such that for any $f\in \Ccinf(\bR^{d+1})$ and $b>0$
\begin{align*}
\aint_{Q_b}\aint_{Q_b}|L f (t,x)-L f(s,y)|dxdtdyds \leq C \|f\|_{L_\infty(\bR^{d+1})}.
\end{align*}
\end{corollary}

\begin{proof}
Take a function $\eta=\eta(t) \in C^\infty(\bR)$ such that $0\leq \eta\leq 1$, $\eta=1$ on $(-\infty, -8\kappa(b)/3)$ and $\eta(t)=0$ for $t\geq -7\kappa(b)/3$. Then for any $(t,x),(s,y)$
\begin{equation*}
\begin{aligned}
|Lf(t,x)-Lf(s,y)| &\leq |Lf(t,x) - Lf(s,x)| + |Lf(s,x) - Lf(s,y)| 
\\
&\leq \sum_{i=1}^{2} \big( |Lf_{i} (t,x) - Lf_{i} (s,x)| + |Lf_{i}(s,x) - Lf_{i}(s,y)| \big),
\end{aligned}
\end{equation*}
where $f_{1}=f\eta$, and $f_{2}=f(1-\eta)$. By Lemma \ref{inwholeestimate} it is easy to see that
\begin{equation*}
\begin{aligned}
&\aint_{Q_{b}}\aint_{Q_{b}} \big( |Lf_{2} (t,x) - Lf_{2} (s,x)| + |Lf_{2}(s,x) - Lf_{2}(s,y)| \big) dxdtdyds 
\\
&\leq 4 \aint_{Q_{b}} |Lf_{2}(t,x)| dxdt
 \leq C \|f\|_{L_{\infty}(\bR^{d+1})}.
\end{aligned}
\end{equation*}
Note that due to Lemma \ref{outwholetimeestimate},
\begin{equation*}
\begin{aligned}
\aint_{Q_{b}}\aint_{Q_{b}} \left|  Lf_{1}(t,x) - Lf_{1}(s,x) \right| dxdtdyds \leq C \|f\|_{L_{\infty}(\bR^{d+1})}.
\end{aligned}
\end{equation*}
Hence, it only remains to consider $|Lf_{1}(s,x) - Lf_{1}(s,y)|$. For this, take $\zeta \in C_c^\infty(\bR^{d})$ such that $0\leq \zeta \leq 1$, $\zeta=1$ on $B_{7b/3}$ and $\zeta=0$ outside of $B_{8b/3}$. Then denoting $f_{3}=f_{1}(1-\zeta)$, and $f_{4}=f_{1}\zeta$
\begin{equation*}
\begin{aligned}
|Lf_{2}(s,x) - Lf_{2}(s,y)| \leq \sum_{i=3}^{4} |Lf_{i}(s,x) - Lf_{i}(s,y)|
\end{aligned}
\end{equation*}
Applying Lemma \ref{outoutspaceestimate} and Lemma \ref{outinspaceestimate} to $f_{3}$ and $f_{4}$ respectively, we have
$$
\aint_{Q_{b}}\aint_{Q_{b}} |Lf_{2}(s,x) - Lf_{2}(s,y)| dxdtdyds \leq C \|f\|_{L_{\infty}(\bR^{d+1})},
$$
and the corollary is proved.
\end{proof}

For  locally integrable functions $f$ on $\bR^{d+1}$, we define the BMO semi-norm of $f$ on $\bR^{d+1}$ as
\begin{equation*}
\|f\|_{BMO(\bR^{d+1})}=\sup_{Q\in \bQ} \aint_Q |f(t,x)-f_Q| dtdx
\end{equation*}
where $f_Q=\aaint_Q f(t,x)dtdx$ and
\begin{equation*}
\bQ:=\{Q_b(t_0,x_0) : b>0, (t_0,x_0)\in \bR^{d+1} \}.
\end{equation*}

For  measurable functions $f$ on $\bR^{d+1}$, we define the sharp function
\begin{align*}
f^{\#}(t,x)=\sup_{Q_{b}(r,z)\in \bQ} \aint_{Q_b(r,z)}|f(s,y)-f_{Q_b(r,z)}|dsdy,
\end{align*}
where the supremum is taken over all $Q_{b}(r,z)\in\bQ$ containing $(t,x)$.

\begin{theorem}[Fefferman-Stein Theorem]
\label{feffermanstein}
For any $1<p<\infty$, and $f\in L_p(\bR^{d+1})$,
\begin{equation*}
C^{-1}\|f^{\#}\|_{L_p(\bR^{d+1})}\leq \|f\|_{L_p(\bR^{d+1})}\leq C\|f^{\#}\|_{L_p(\bR^{d+1})},
\end{equation*}
where $C>1$ depends on $\alpha,d,p,\kappa_{1},\kappa_{2}$.
\end{theorem}

\begin{proof}
See  \cite[Theorem I.3.1, Theorem IV.2.2]{stein1993harmonic}.  We only remark that  due to \eqref{eqn 05.12:15:07},  the balls $Q_b(s,y)$ satisfy the conditions (i)--(iv) in  \cite[Section 1.1] {stein1993harmonic}.
\end{proof}

Recall that linear operators $L$  is given by
$$
L f(t,x)=\lim_{\varepsilon \downarrow 0} \int_{-\infty}^{t-\varepsilon} \left(\int_{\bR^d} q^{\alpha,\alpha+1}_{d}(t-s,x-y) f(s,y) dy\right) ds,
$$
The following theorem is our main result in this section. The proof is quite standard.

\begin{theorem}\label{thm 05.14.18:28}

(i) For any $f\in L_2(\bR^{d+1})\cap L_\infty(\bR^{d+1})$,
\begin{equation*}
\begin{gathered}
\|L f\|_{BMO(\bR^{d+1})} \leq C(\alpha,d,\kappa_{1},\kappa_{2},\text{\boldmath{$\delta$}},\ell,C_{0},C_{1},C_{2}) \|f\|_{L_\infty(\bR^{d+1})}.
\end{gathered}
\end{equation*}

(ii) For any $p,q\in(1,\infty)$ and $f\in C_c^\infty (\bR^{d+1})$,
\begin{equation*}
\begin{gathered}
\|L f\|_{L_q(\bR;L_p(\R^{d}))} \leq C(\alpha,d,p,q,\kappa_{1},\kappa_{2},\text{\boldmath{$\delta$}},\ell,C_{0},C_{1},C_{2}) \|f\|_{L_q(\bR;L_p(\R^{d}))}.
\end{gathered}
\end{equation*}

\end{theorem}

\begin{proof}
If $\mathcal{L}_{X}=-\phi(-\Delta)$ for Bernstein function $\phi$ satisfying \eqref{phi_lsc}, then the theorem is a direct consequence of \cite[Theorem 4.10]{kim2020nonlocal}. Since we have all the necessary results (i.e. Corollary \ref{qintegrable} and Corollary \ref{outwholeestimate}) we can prove theorem by following the proof of \cite[Theorem 4.10]{kim2020nonlocal}.
\end{proof}

\medskip

\mysection{Proof of Theorem \ref{main theorem2}}\label{s:thm_proof}

In this section, we will prove Theorem \ref{main theorem2}.

\vspace{1mm}
\noindent\textbf{Proof of Theorem \ref{main theorem2}.} 
Due to Lemma \ref{basicproperty}, it suffices to prove  case $\gamma=0$.

\vspace{1mm}

\noindent{\textbf{Step 1}} (Existence and estimation of solution).

First assume $f\in C_c^\infty(\bR^{d+1}_+)$, and let $u(t,x)$ be a function with representation \eqref{u=qf}. Using Remark \ref{Hvaluedconti} and the integrability of $q^{\alpha,1}_{d}$, we can easily check 
 $D^m_x u$, $\mathcal{L}D^m_x u$ $\in C([0,T];L_p)$, and thus $u\in C^{\alpha,\infty}_p([0,T]\times\bR^d)$. Also, by Lemma \ref{u=qfsolution},  $u$ satisfies equation \eqref{mainequation} with $u(0,\cdot)=0$. 
 
Now we show estimation \eqref{mainestimate2} and \eqref{eqn 05.27.14:12}. Take $\eta_{k}=\eta_{k}(t)\in C^{\infty}(\bR)$ such that $0\leq \eta_{k}\leq 1$, $\eta_{k}(t)=1$ for $t\leq T+1/k$ and $\eta_{k}(t)=0$ for $t\geq T+2/k$. Since $f\eta_{k}\in L_{q}(\bR;L_{p}(\bR^{d}))$, and $f(t)=f\eta_{k}(t)$ for $t\leq T$, By Theorem \ref{thm 05.14.18:28} (ii), we have
\begin{equation*}
\begin{aligned}
\|\mathcal{L}u\|_{\mathbb{L}_{q,p}(T)} &= \|L f \|_{\mathbb{L}_{q,p}(T)} = \|L(f\eta_{k})\|_{\mathbb{L}_{q,p}(T)} 
\\
&\leq \|L(f\eta_{k})\|_{L_{q}(\bR;L_{p}(\bR^{d}))}\leq C \|f\eta_{k}\|_{L_{q}(\bR;L_{p}(\bR^{d}))}.
\end{aligned}
\end{equation*}
Hence, by the dominated convergence theorem, taking $k\to\infty$, we have
$$
\|\mathcal{L}u\|_{\mathbb{L}_{q,p}(T)} \leq C \|f\|_{\mathbb{L}_{q,p}(T)}.
$$
Also, by Corollary \ref{qintegrable} and Minkowski's inequality, we can easily check that
$$
\|u\|_{\mathbb{L}_{q,p}(T)} \leq C(T) \|f\|_{\mathbb{L}_{q,p}(T)}.
$$
Therefore, using the above inequalities and \eqref{eqn 03.25.15:03} we prove estimation  \eqref{mainestimate2} and \eqref{eqn 05.27.14:12}.
For general $f$,  we take a sequence of functions   $f_{n}\in \Ccinf(\R^{d+1}_{+})$  such that $f_n \to f$ in $\bL_{q,p}(T)$. Let $u_{n}$  denote the solution with representation  \eqref{u=qf} with $f_{n}$ in place of $f$. Then \eqref{mainestimate2} applied to $u_m-u_n$ shows that $u_{n}$ is a Cauchy sequence in 
$\bH^{\alpha,\psi,2}_{q,p,0}(T)$. By taking $u$ as the limit of $u_{n}$ in $\bH^{\alpha,\psi,2}_{q,p,0}(T)$, we find that $u$ satisfies \eqref{mainequation2} and estimation  \eqref{mainestimate2} and \eqref{eqn 05.27.14:12} also holds for $u$. 

\vspace{1mm}

\noindent\textbf{Step 2} (Uniqueness  of solution).

Let  $u\in \bH_{q,p,0}^{\alpha,\psi,2}(T)$ be a solution to equation \eqref{mainequation2} with $f=0$. Take $u_{n}\in \Ccinf(\bR^{d+1}_{+})$ which converges to $u$ in $\mathbb{H}^{\alpha,\psi,2}_{q,p}(T)$, and let $f_{n}:=\partial^{\alpha}_{t}u_{n}-\mathcal{L}u_{n}$. Then by Lemma \ref{u=qfsolution}, $u_{n}$ satisfies representation \eqref{u=qf} with $f_{n}$. Hence, by the argument in \textbf{Step 1}, we have
$$
\|u_{n}\|_{\mathbb{H}^{\psi,2}_{q,p}(T)} \leq C(T) \leq \|f_{n}\|_{\mathbb{L}_{q,p}(T)}.
$$
Since $f_{n}=\partial^{\alpha}_{t}u_{n}-\mathcal{L}u_{n}$ converges to $0$ due to the choice of $u_{n}$, we conclude that $u=0$ in $\mathbb{H}^{\alpha,\psi,2}_{q,p}(T)$. The theorem is proved.

\qed

\mysection{Further discussions}\label{sec 11.01.14:33}
The purpose of this section is to handle the case $\ell=\ell_{d}$ depends on $d$. A typical example is $\mathcal{L}=-\log{(1-\Delta)}$. In this case, the corresponding jumping kernel $j_{d}$ satisfies
$$
j_{d}(r) \asymp r^{-d}\ell_{d}(r^{-1}) \quad \text{for}\;\, r>0,
$$
where 
\begin{equation}\label{eqn 11.02.15:18}
\ell_{d}(r) := r^{(1-d)/2}e^{-r^{-1}} + e^{-1}\mathbf{1}_{r\geq 1}
\end{equation}
(see \cite{SSV06}).  To handle operators whose jump intensity $\ell_{d}$ depend on $d$, we impose the following.

\begin{assumption}\label{asm 11.01.15:05}
(i) The function $-\frac{1}{r}\left(\frac{d}{dr}j_{d}\right)(r)$ is decreasing and there exist continuous functions $\ell_{d}$ and $\ell_{d+2}$ such that
\begin{equation}\label{eqn 11.02.15:10}
\kappa_{1}r^{-d}\ell_{d}(r^{-1}) \leq j_{d}(r) \leq \kappa_{2} r^{-d}\ell_{d}(r^{-1}) \quad \text{for} \quad r>0,
\end{equation}
and
\begin{equation}\label{eqn 11.02.15:12}
\kappa_{1} r^{-d-2} \ell_{d+2}(r^{-1}) \leq -\frac{1}{r}\frac{d}{dr} j_{d}(r) \leq \kappa_{2} r^{-d-2}\ell_{d+2}(r^{-1}) \quad \text{for} \quad r>0.
\end{equation}
Moreover, $\ell_{d}$ and $\ell_{d+2}$ satisfy \eqref{H:s} and $\limsup_{r\to\infty}(\ell_{d}(r) \vee \ell_{d+2}(r))<\infty$.

(ii) $h_{d}$ satisfies \eqref{H:l} with $\delta_{3}>0$ and
$$
K_{d}(r)  \asymp K_{d+2}(r), \quad  h_{d}(r) \asymp h_{d+2}(r) \quad \text{for} \quad r>0,
$$
where $K_{d+i}$ and $h_{d+i}$  are defines as in \eqref{eqn 11.02.14:53} with $\ell_{d+i}$ ($i=0,2$).
\end{assumption}

\begin{remark}
(i) When the given process $X_{t}$ is a subordinate Brownian motion, then \eqref{eqn 11.02.15:10} implies \eqref{eqn 11.02.15:12}. See Remark \ref{r:assum} (v).

(ii) It is easy to see that $\ell_{d}(r) = r^{(1-d)/2}e^{-r^{-1}} + e^{-1}\mathbf{1}_{r\geq 1}$ from \eqref{eqn 11.02.15:18} satisfies \eqref{H:s} with $\delta_{1}=\delta_{2}=0$. Also, we can check that (recall \eqref{eqn 05.27.15:40})
$$
K_{d}(r) \asymp (1\wedge r^{-2}) \asymp K_{d+2}(r),  \quad  h_{d}(r) \asymp \log{(1+r^{-2})} \asymp h_{d+2}(r) \quad \text{for} \quad r>0.
$$
In this case, the corresponding process $X_{t}$ is a subordinate Brownian motion since $\phi(r)=\log{(1+r)}$ is a Bernstein function. Hence,  Assumption \ref{asm 11.01.15:05} holds with $\ell_{d}(r) = r^{(1-d)/2}e^{-r^{-1}} + e^{-1}\mathbf{1}_{r\geq 1}$.

(iii) If $\limsup_{r\to\infty}\ell(r)=\infty$, then by Proposition \ref{p:hku} (ii) we have
$$
p_{d}(t,x) \leq c_{2}t\frac{K(\theta_{a_0}(|x|,t))}{[\theta_{a_0}(|x|,t)]^d}\exp{(-b_{2}th(\theta_{a_0}(|x|,t)))}.
$$
In this case $\theta_{a_{0}}(r,t) = r \vee (\ell_{d}^{-1}(a_{0}/t))^{-1}$ may depend on $d$, and we need additional assumption to follow our approach.
\end{remark}

\begin{theorem}
Let $\alpha\in(0,1)$, $p,q\in(1,\infty)$, $q'=q/(q-1)$, $\gamma \in \bR$, and $T\in(0,\infty)$. Suppose Assumption \ref{asm 11.01.15:05} holds. Then for any  $f\in \bH_{q,p}^{\psi,\gamma}(T)$, the equation
\begin{equation*}
\partial_t^\alpha u = \mathcal{L}u + f,\quad t>0\,; \quad u(0,\cdot)=0
\end{equation*}
has a unique solution $u$ in the class $\bH_{q,p,0}^{\alpha,\psi,\gamma+2}(T)$, and for the solution $u$ it holds that
\begin{equation*}
\|u\|_{\bH_{q,p}^{\alpha,\psi,\gamma+2}(T)}\leq C \|f\|_{\bH_{q,p}^{\psi,\gamma}(T)},
\end{equation*}
where $C>0$ depends only on $\alpha,d,\kappa_{1},\kappa_{2},p,q,\gamma,\ell,C_{0},C_{1},C_{2},T$, and {\boldmath$\delta$}. Furthermore, we have
\begin{equation*}
\| \mathcal{L} u\|_{\mathbb{H}^{\psi,\gamma}_{q,p}(T)} \leq C \| f\|_{\mathbb{H}^{\psi,\gamma}_{q,p}(T)},
\end{equation*}
where $C>0$ depends only on $\alpha,d,\kappa_{1},\kappa_{2},p,q,\gamma,\ell,C_{0},C_{1},C_{2}$ and {\boldmath$\delta$}.
%
%(ii) Suppose that $\mathcal{L}=-\phi(-\Delta)$, where $\phi$ is a Bernstein function. Suppose further that Assumption \ref{asm 11.01.15:05} holds and that $h_{d}$ satisfies Assumption \ref{asm 11.07.16:00}. Then for any  $f\in \bH_{q,p}^{\psi,\gamma}(T)$ and $u_{0}\in B^{\psi,\gamma+2+2\varepsilon_{0}/q'-2/\alpha q}_{p,q}$, the equation
%\begin{equation*}
%\partial_t^\alpha u = \mathcal{L}u + f,\quad t>0\,; \quad u(0,\cdot)=u_{0}
%\end{equation*}
%has a unique solution $u$ in the class $\bH_{q,p}^{\alpha,\psi,\gamma+2}(T)$, and for the solution $u$ it holds that
%\begin{equation*}
%\|u\|_{\bH_{q,p}^{\alpha,\psi,\gamma+2}(T)}\leq C  \left( \|u_{0}\|_{B^{\psi,\gamma+2+2\varepsilon_{0}/q'-2/\alpha q}_{p,q}}  +  \|f\|_{\bH_{q,p}^{\psi,\gamma}(T)}  \right),
%\end{equation*}
%where $C>0$ depends only on $\alpha,d,\kappa_{1},\kappa_{2},p,q,\gamma,\varepsilon_{0},\ell,C_{0},C_{1},C_{2},T$ and {\boldmath$\delta$}.
\end{theorem}
\begin{proof}
Under Assumption \ref{asm 11.01.15:05} (i) we have Proposition \ref{p:hku} (i) and Lemma \ref{l:hke_largetime}. Also, by Assumption \ref{asm 11.01.15:05} (ii) we may assume that $K$ and $h$ are independent of $d$. Hence, by following the argument in Section \ref{sec 4}, Section \ref{sec 5} and Section \ref{s:thm_proof} we obtain the theorem.
\end{proof}

\appendix

\mysection{Appendix}\label{s:App}
\begin{lemma}\label{lem 09.01.17:25}
(i) Suppsose that $\ell_{1}\in \cG$ and let $\ell_{2}$ be a function which satisfies $\ell_{2} \in (a,b)$ for some $0<a<b<\infty$. Then $\ell_{1}\ell_{2}\in \cG$.

(ii) Suppose $\ell\in\cG$. Then for any $b>0$, $r\mapsto\ell(r^b)$ belongs to $\cG$.

(iii) Let $\ell_1\in\cG$ and $\ell_2$ be an increasing function satisfying $\ell_2\ge c$ on $[1,\infty)$ for some $c>0$. Then, $\ell_1/\ell_2\in\cG$.

(iv) Let $\ell_{1}(r)=\log{(1+r)}$ and let $\ell_{k+1}=\ell_{k}\circ\ell_{1}(r)$ for $k\in\bN$. Then for any $n\in \bN$ $k_{1},\dots,k_{n}\in \bN$ and $b_{1},\dots,b_{n}>0$ we have
$$
\Lambda (r) = \prod_{i=1}^{n} (\ell_{k_{i}}(r))^{b_{i}} \in \cG.
$$

(v) For $b\in(0,1/2)$, $\ell(r)=(e^{(\log{(1+r)})^{b}}-1)  \in \cG$.
\end{lemma}
\begin{proof}
(i)$\&$(ii) Trivial.

(iii) Let $\ell=\ell_1/\ell_2$. Using that $\ell_2\ge c$, $\ell_2$ is increasing and $\ell_1\in\cG$, we see that for any $a>0$
\begin{align*}
&\sup_{r>1} \int^{r}_{1}\frac{\ell_1(s)}{s\ell_2(s)}ds\cdot\exp{\left( -\frac{a\ell_2(r)}{\ell_1(r)}\int^{r}_{1}\frac{\ell_1(s)}{s\ell_2(s)}ds  \right)}\\
&\le\sup_{r>1}c^{-1} \int^{r}_{1}\frac{\ell_1(s)}{s}ds\cdot\exp{\left( -\frac{a}{\ell_1(r)}\int^{r}_{1}\frac{\ell_1(s)}{s}ds  \right)} \leq C.
\end{align*}

(iv) Let $\tilde{\Lambda}=(\ell_{1})^{2}\Lambda$. Then $\tilde{\Lambda}(r) = \prod^{n+1}_{i=1} (\ell_{k_{i}}(r))^{b_{i}}$,
and we may set $k_{1}=1$, $b_{1}=2$. If we show that $\tilde{\Lambda} \in \cG$, then due to (iii), it follows that $\Lambda\in \cG$. Thus we will show that $\tilde{\Lambda} \in \cG$. Set $\ell_{0}(s)=s$. By using the change of variable, 
\begin{align*}
\int_{1}^{r}\tilde{\Lambda}(s) \, s^{-1} ds &\asymp \int_{\log2}^{\log{(1+r)}} \prod_{i=1}^{n+1} (\ell_{k_{i}-1}(s))^{b_{i}} ds.
\end{align*}
It is easy to see that
\begin{align*}
\int_{\log2}^{\log{(1+r)}} \prod_{i=1}^{n+1} (\ell_{k_{i}-1}(s))^{b_{i}} ds \leq C \log{(1+r)}\tilde{\Lambda}(r).
\end{align*}
Moreover, by the integration by parts and  $((\ell_{k_{i}-1})^{b_{i}})'(s) \leq C(b_{i},k_{i})$,
\begin{align*}
&\int_{\log2}^{\log{(1+r)}} \prod_{i=1}^{n+1} (\ell_{k_{i}-1}(s))^{b_{i}} ds
\\
& = \log{(1+r)}\tilde{\Lambda}(r)-\log{2}\tilde{\Lambda}(1) - \int_{\log2}^{\log{(1+r)}} \left( \prod_{i=1}^{n+1} (\ell_{k_{i}-1}(\cdot))^{b_{i}} \right)'(s) \, s\,ds
\\
& \geq \log{(1+r)}\tilde{\Lambda}(r)-\log{2}\tilde{\Lambda}(1) - C \int_{\log2}^{\log{(1+r)}} s\, ds
\\
& \geq \log{(1+r)}\tilde{\Lambda}(r)-\log{2}\tilde{\Lambda}(1) - C (\log{(1+r)})^{2}.
\end{align*}
Since $k_{1}=1$ and $b_{1}=2$ (for $\tilde{\Lambda}$),  there exists $M>1$ such that
$$\int_{1}^{r}\tilde{\Lambda}(s)\, s^{-1} ds \geq \frac{1}{2} \log{(1+r)}\tilde{\Lambda}(r),\quad\text{for}\;\;r>M.$$
Thus, we have
\begin{align*}
& \sup_{r>M} \int_{1}^{r} \tilde{\Lambda}(s) \,s \,ds \cdot \exp{\left(-\frac{a}{\tilde{\Lambda}(r)} \int_{1}^{r} \tilde{\Lambda}(s) \,s \,ds \right)}
\\
&\leq C \sup_{r>M}  \log{(1+r)}\tilde{\Lambda}(r) \cdot \exp{\left( -c \log{(1+r)}   \right)}
 \leq C.
\end{align*}
The above inequality implies that $\Lambda\in\cG$ since it is bounded on $(1, M]$. 

(v) By the change of variable, we have
$$
\int_{1}^{r} \ell(s) s^{-1} \asymp \int_{\log2}^{\log{(1+r)}} (e^{s^{b}}-1) ds  \leq C e^{\left( (\log{(1+r)})^{b} \right)} \log{(1+r)}.
$$
Also, by the integration by parts, we obtain
\begin{align*}
\int_{\log2}^{\log{(1+r)}} (e^{s^{b}}  - 1) ds 
&= \ell(r) \log{(1+r)} - \ell(1)\log2 -  b  \int_{\log2}^{\log{(1+r)}} s^{b} e^{s^{b}} ds
\\
& \geq \ell(r) \log{(1+r)} - \ell(1)\log2 
\\
&\quad- b (\log{(1+r)})^{b}  \int_{\log2}^{\log{(1+r)}} (e^{s^{b}}-1) ds - b(\log{(1+r)})^{b+1} .
\end{align*}
This shows that there exists $M>1$ such that for $r>M$
$$
\frac{1}{\ell(r)} \int_{\log2}^{\log{(1+r)}} (e^{s^{b}}  - 1) ds \geq  (\log{(1+r)})^{1-b} -C .
$$
Hence, we have (recall that $b\in(0,1/2))$
\begin{align*}
&\sup_{r>M} \int^{r}_{1}\frac{\ell(s)}{s}ds\cdot\exp{\left( -\frac{a}{\ell(r)}\int^{r}_{1}\frac{\ell(s)}{s}ds  \right)}
\\
& \leq C \sup_{r>M} \left( e^{\left( (\log{(1+r)})^{b} \right)} \log{(1+r)} e^{ \left(- c (\log{(1+r)})^{1-b} \right)} \right)
\leq C.
\end{align*}
Thus $\ell\in \cG$. The lemma is proved.
\end{proof}

\begin{lemma}\label{lem 07.23.16:47}
Let $f:(0,\infty) \to (0,\infty)$ be an increasing continuous function satisfying $\sup_{0<r<1}f(r)<\infty$, $\lim_{r\to\infty}f(r)=\infty$ and
\begin{equation*}
 \frac{f(R)}{f(r)} \leq c_{1} \left( \frac{R}{r} \right)^{\delta} \quad \text{for}\;\; 1\leq r\leq R<\infty
\end{equation*}
for some $c_{1}, \delta>0$. Then there is a strictly increasing continuous function $\tilde{f}:(0,\infty) \to (0,\infty)$ satisfying
$$
f(r) \leq \tilde{f}(r) \leq C f(r) \quad \text{for}\; \, r>0,
$$ 
where the constant $C$ does not depend on $r$.
\end{lemma}
\begin{proof}
We prove the lemma by constructing $\tilde{f}$. Extend $f$ to $\bR_{+} \cup \{0\}$ by 
$$
f(0) = \lim_{r\downarrow 0 } f(r).
$$
Now let $A= \{ r\geq0 : \exists\, s\geq0, s\neq r ~ \text{such that} ~ f(r)=f(s) \}$.
Then, we can check that $A=\cup_{k=1}^{\infty}[r_{k},l_{k}]$, where $[r_{k},l_{k}]$ are pairwise disjoint closed intervals.

\textbf{Case 1.} Assume that $f(0)=0$ or $f(r)$ is strictly increasing for $r\leq 1$. Then, there is a positive number $a>0$ such that $A\subset [0,a)^{c}$. Note that for $[r_{k},l_{k}]$ we can choose $\varepsilon_{k}\in(0,1)$ such that $(l_{k},l_{k}+\varepsilon_{k}]\subset A^{c}$. Now define $\tilde{f}_{k}$ on $[r_{k},r_{l}+\varepsilon_{k}]$ as
\begin{equation*}
\tilde{f}_{k}(r) =  \frac{f(l_{k}+\varepsilon_{k})-f(r_{k})}{l_{k}+\varepsilon_{k}-r_{k}}(r-r_{k}) + f(r_{k}).
\end{equation*}
Then $\tilde{f}_{k}$ is continuous, strictly increasing on $[r_{k},l_{k}+\varepsilon_{k}]$ and it satisfies $\tilde{f}_{k}(r_{k})=f(r_{k})$, and $\tilde{f}_{k}(l_{k}+\varepsilon_{k}) = f(l_{k}+\varepsilon_{k})$. Moreover, on $[r_{k},l_{k}+\varepsilon_{k}]$, $\tilde{f}_{k}$ satisfies
$$
1 \leq \frac{\tilde{f}_{k}(r)}{f(r)} \leq \frac{f(l_{k}+\varepsilon_{k})}{f(r_{k})} \leq c_{1} \left ( \frac{l_{k}+\varepsilon_{k}}{l_{k}} \right)^{\delta} \leq C(a)
$$
since $f(r_{k})=f(l_{k})$, $l_{k}>r_{k}\geq a$. Now define $\tilde{A} = \cup_{k=1}^{\infty}[r_{k},l_{k}+\varepsilon_{k}]$ and
\begin{equation}\label{eqn 08.18.10:25}
\tilde{f}(r) = \mathbf{1}_{\tilde{A}^{c}}(r)f(r) +  \sum_{k=1}^{\infty} \mathbf{1}_{[r_{k},l_{k}+\varepsilon_{k}]} \tilde{f}_{k}(r).
\end{equation}
Then $\tilde{f}$ is a desired function.

\textbf{Case 2.} Now assume that $f(0)\neq 0$ and $f(r)$ is not strictly increasing for $r\leq 1$. Then there exists $b\geq 1$ such that $[0,b]=[r_{1},l_{1}]$. Take $\varepsilon_{1}$ and $\tilde{f}_{1}$ for $[r_{1},l_{1}]$ corresponding to $\varepsilon_{k}$ and $\tilde{f}_{k}$ in above case. Then on $[r_{1},l_{1}]$, we have
$$
1 \leq \frac{\tilde{f}_{1}(r)}{f(r)} \leq \frac{f(b+\varepsilon_{1})}{f(0)} \leq c_{1} \left ( \frac{b+\varepsilon}{b} \right)^{\delta} \leq C(b)
$$
since $f(0)=f(b)>0$. For other $k$, we have the same result by following \textbf{Case 1}. Hence, by taking $\tilde{f}$ as in \eqref{eqn 08.18.10:25}, the lemma is proved.
\end{proof}

\begin{lemma}\label{l:qest_offd}
Let  $f:(0,\infty)\to(0,\infty)$ be a strictly increasing continuous function and $f^{-1}$ be its inverse. Suppose that there exist $c,\gamma>0$ such that $(f(R)/f(r))\le c(R/r)^\gamma$ for $0<r\le R<\infty$. Then, for any $k>0$, there exists $C>0$ such that for any $b>0$
\begin{align*}
\int_{(f(b^{-1}))^{-1}}^\infty s^{-1}f^{-1}(s^{-1})^k ds \le Cb^{-k}.
\end{align*}
\end{lemma}

\begin{proof}
By the scaling property of $f$ with $R=b^{-1}$ and $r=f^{-1}(s^{-1})$, and the fact that $f(f^{-1}(s^{-1}))= s^{-1}$, we see that
\begin{align*}
\frac{f(b^{-1})}{s^{-1}}  \le c\left(\frac{b^{-1}}{f^{-1}(s^{-1})}\right)^{\gamma} \quad \text{for} \quad s>(f(b^{-1}))^{-1}.
\end{align*}
Thus,
\begin{equation*}
\begin{aligned}
\int_{(f(b^{-1}))^{-1}}^\infty s^{-1}f^{-1}(s^{-1})^k ds
&\leq Cb^{-k}f(b^{-1})^{-k/\gamma} \int_{(f(b^{-1}))^{-1}}^\infty s^{-1-k/\gamma}    ds = Cb^{-k}.
\end{aligned}
\end{equation*}
\end{proof}

We use the following lemma with $f(r) = h(r^{-1})$. Note that by using \eqref{eqn 05.12:15:07}, one can check that $h(r^{-1})$ is a strictly increasing function satisfying $h(R^{-1})\leq c (R/r)^{2} h(r^{-1})$ for any $0<r<R$.
\begin{lemma}\label{l:qest_int}
Let  $f:(0,\infty)\to(0,\infty)$ be a strictly increasing function and $f^{-1}$ be its inverse. Suppose that there exist $c,\gamma>0$ such that $(f(R)/f(r))\le c(R/r)^\gamma$ for $0<r\le R<\infty$. Then, there exists $C>0$ such that for any $b>0$
\begin{align*}
\int_{(f(b^{-1}))^{-1/\alpha}}^\infty \int_{b\leq |y| \leq (f^{-1}(s^{-\alpha}))^{-1}} \int_{(f(|y|^{-1}))^{-1}}^{2s^{\alpha}} (f^{-1}(r^{-1}))^{d+1}  s^{-\alpha-1} dr dy ds\le Cb^{-1}.
\end{align*}
\end{lemma}

\begin{proof}
By the change of variable and Fubini's theorem and Lemma \ref{l:qest_offd},
\begin{align*}
&\int_{(f(b^{-1}))^{-1/\alpha}}^\infty \int_{b\leq |y| \leq (f^{-1}(s^{-\alpha}))^{-1}} \int_{(f(|y|^{-1}))^{-1}}^{2s^{\alpha}} f^{-1}(r^{-1})^{d+1}  s^{-\alpha-1} dr dy ds\\
&=C\int_{(f(b^{-1}))^{-1/\alpha}}^\infty \int_{b}^{(f^{-1}(s^{-\alpha}))^{-1}} \int_{(f(\rho^{-1}))^{-1}}^{2s^{\alpha}} f^{-1}(r^{-1})^{d+1}  s^{-\alpha-1}\rho^{d-1} dr d\rho ds\\
&\le C\int_{(f(b^{-1}))^{-1/\alpha}}^\infty \int_{(f(b^{-1}))^{-1}}^{2s^{\alpha}}\int_{0}^{(f^{-1}(r^{-1}))^{-1}}  f^{-1}(r^{-1})^{d+1}  s^{-\alpha-1}\rho^{d-1} d\rho dr  ds\\
&\le C\int_{(f(b^{-1}))^{-1/\alpha}}^\infty \int_{(f(b^{-1}))^{-1}}^{2s^{\alpha}} f^{-1}(r^{-1})  s^{-\alpha-1} dr  ds\\
&\le C\int_{(f(b^{-1}))^{-1}}^\infty \int_{(r/2)^{1/\alpha}}^{\infty} f^{-1}(r^{-1})  s^{-\alpha-1}  ds dr \\
&\le C\int_{(f(b^{-1}))^{-1}}^\infty r^{-1}f^{-1}(r^{-1}) dr \le Cb^{-1}.
\end{align*}
\end{proof}

\begin{lemma}\label{l:qest_case1_pf} 
Let $\alpha\in(0,1)$. Suppose the function $\ell$ satisfies Assumption \ref{ell_con} (i) let $\kappa(b) = (h(b))^{-1/\alpha}$ and let $t_{1}>0$ be taken from Lemma \ref{l:hke_largetime}. Then,
there exists $C>0$ depending only on $\alpha,\kappa_{1},\kappa_{2},d,\ell,C_{0},C_{1},C_{2}$, and {\boldmath{$\delta$}} such that for any $b>0$
\begin{align}
&\int_{\kappa(b)}^\infty \int_{b\leq |y| \leq h^{-1}(s^{-\alpha})} \int_{0}^{t_1} |D_{x}p_{d} (r,y)| |\varphi_{\alpha,\alpha+1} (s,r)| dr dy ds\le Cb^{-1},\label{q1upper_case1_pf}\\
&\int_{\kappa(4b)}^{\infty} \int_{|y|\leq 4b} \int_{0}^{t_1} p_{d} (r,y)|\varphi_{\alpha,\alpha+1} (s,r)| dr dy ds\le C\label{q0upper_case1_pf}.
\end{align}
\end{lemma}

\begin{proof}
By \eqref{p'_est}, \eqref{j'_est}, Proposition \ref{p:hku} (i), \eqref{philarge} and \eqref{betainteger},
\begin{equation*}
\begin{aligned}
&\int_{\kappa(b)}^\infty \int_{b\leq |y| \leq h^{-1}(s^{-\alpha})} \int_{0}^{t_1} |D_{x}p_{d} (r,y)| |\varphi_{\alpha,\alpha+1} (s,r)| dr dy ds
\\
&\leq \int_{\kappa(b)}^\infty \int_{b}^{h^{-1}(s^{-\alpha})} \int_{0}^{\infty} rK(\rho) \rho^{-2} e^{-Crh(\rho)}s^{-\alpha-1} dr d\rho  ds
\\
&\leq C \int_{(h(b))^{-1/\alpha}}^\infty \int_{b}^{h^{-1}(s^{-\alpha})} \rho^{-2} \frac{1}{h(\rho)} s^{-\alpha-1} d\rho ds,
\end{aligned}
\end{equation*}
where we used the relations $se^{-s} \leq Ce^{-s/2}$ ($s>0$) and $K(\rho)\leq h(\rho)$ for the last inequality. By Fubini's theorem we have
\begin{equation}\label{eqn 05.30.18:42}
\begin{aligned}
&\int_{\kappa(b)}^\infty \int_{b\leq |y| \leq h^{-1}(s^{-\alpha})}  \int_{0}^{t_1} |D_{x}p_{d} (r,y)| |\varphi_{\alpha,\alpha+1} (s,r)| dr dy ds
\\
&\leq C \int_{(h(b))^{-1/\alpha}}^\infty  \int_{b}^{h^{-1}(s^{-\alpha})} \rho^{-2} \frac{1}{h(\rho)} s^{-\alpha-1} d\rho ds
\\
&\leq C \int_{b}^{\infty}\int_{(h(\rho))^{-1/\alpha}}^{\infty} s^{-\alpha-1} \frac{1}{h(\rho)} \rho^{-2} ds d\rho
\leq C \int_{b}^{\infty} \rho^{-2} d\rho
= C b^{-1},
\end{aligned}
\end{equation}
which shows \eqref{q1upper_case1_pf}.

Now we prove \eqref{q0upper_case1_pf}. Using Proposition \ref{p:hku} (i), \eqref{philarge} and \eqref{betainteger}, we see that
\begin{equation}\label{eqn 08.05.11:51}
\begin{aligned}
&\int_{\kappa(4b)}^{\infty} \int_{|y|\leq 4b} \int_{0}^{t_1} p_{d} (r,y)\varphi_{\alpha,\alpha+1} (s,r) dr dy ds
\\
&\leq C \int_{\kappa(4b)}^{\infty} \int_{0}^{4b} \int_{0}^{\infty} rK(\rho) \rho^{-1} e^{-C^{-1}rh(\rho)} s^{-\alpha-1}dr d\rho ds
\\
&\leq C \int_{\kappa(4b)}^{\infty}  \int_{0}^{\infty} \int_{\bR^{d}} rK(|y|) |y|^{-d} e^{-C^{-1}rh(|y|)} e^{-C^{-1}rh(4b)}  s^{-\alpha-1} dy dr  ds
\\
&\leq C \int_{\kappa(4b)}^{\infty}  \int_{0}^{\infty} e^{-C^{-1}rh(4b)/2} s^{-\alpha-1}dr ds
\\
&\leq C \int_{(h(4b))^{-1/\alpha}}^{\infty} (h(4b))^{-1} s^{-\alpha-1} ds 
\leq C,
\end{aligned}
\end{equation}
where for the third inequality we use Lemma \ref{rmk 05.07.16:12} (i).
\end{proof}

\medskip

The following lemma is counter part of Lemma \ref{l:qest_case1_pf}. The proof is more delicate than that of Lemma \ref{l:qest_case1_pf} due to the fact that $h(r)$ and $\ell(r^{-1})$ may not be comparable for $0<r\le1$. 

\begin{lemma} 
Let $\alpha\in(0,1)$. Suppose the function $\ell$ satisfies Assumption \ref{ell_con} (ii)--(2). Let $\kappa(b)=(h(b))^{-1/\alpha}$ and let $t_{1}>0$ be taken from Lemma \ref{l:hke_largetime}. Then,
there exists $C>0$ depending only on $\alpha,\kappa_{1},\kappa_{2},d,\ell,C_{0},C_{1},C_{2}$, and {\boldmath{$\delta$}} such that for any $b>0$
\begin{align}
&\int_{\kappa(b)}^\infty \int_{ b\leq |y| \leq h^{-1}(s^{-\alpha})} \int_{0}^{t_1} |D_{x}p_{d} (r,y)| |\varphi_{\alpha,\alpha+1} (s,r)| dr dy ds\le Cb^{-1},\label{q1upper_case2_pf}\\
&\int_{\kappa(4b)}^{\infty} \int_{|y|\leq 4b} \int_{0}^{t_1} p_{d} (r,y)|\varphi_{\alpha,\alpha+1} (s,r)| dr dy ds\le C.\label{q0upper_case2_pf}
\end{align}
\end{lemma}

\begin{proof}
Note that $p(t,0)$ is well-defined on $(0,t_1]$. We first show \eqref{q0upper_case2_pf}.
We split the integral into two parts.
\begin{align*}
\begin{split}
&\int_{\kappa(4b)}^{\infty} \int_{|y|\leq 4b}\int_{0}^{t_1}p_{d} (r,y) |\varphi_{\alpha,\alpha+1} (s,r)| dr dy ds\\
&\leq \int_{\kappa(4b)}^{\infty} \int_{|y|\leq 4b}\int_{0}^{a_0(\ell^*(|y|^{-1}))^{-1}}p_{d}(r,y) |\varphi_{\alpha,\alpha+1} (s,r)| dr dy ds\\
&\quad+\int_{\kappa(4b)}^{\infty} \int_{|y|\leq 4b}\mathbf{1}_{a_0(\ell^*(|y|^{-1}))^{-1} \leq t_{1}} \int_{a_0(\ell^*(|y|^{-1}))^{-1}}^{t_1}p_{d}(r,y)  |\varphi_{\alpha,\alpha+1} (s,r)| dr dy ds\\
&=:I+II.
\end{split}
\end{align*}
We can obtain $I\le C$ by using Proposition \ref{p:hku} (ii) and the same argument in the proof of \eqref{q0upper_case1_pf} (see \eqref{eqn 08.05.11:51}). Thus, we will show $II\le C$ for some constant $C$ for the rest of the proof of \eqref{q0upper_case2_pf}. 
Observe that
\begin{align*}
II&\le\int_{\kappa(4b)}^{\infty} \int_{|y|\leq 4b}\int_{a_0(\ell^*(|y|^{-1}))^{-1}}^{t_1}p_{d}(r,0)  |\varphi_{\alpha,\alpha+1} (s,r)| dr dy ds
\\
&\leq \int_{\kappa(4b)}^{\infty} \int_{|y|\leq 4b}\int_{a_0(\ell^*(|y|^{-1}))^{-1}}^{a_0(\ell^*((4b)^{-1}))^{-1}}p_{d}(r,0)  |\varphi_{\alpha,\alpha+1} (s,r)| dr dy ds
\\
& \quad + \int_{\kappa(4b)}^{\infty} \int_{|y|\leq 4b} \mathbf{1}_{a_0(\ell^*((4b)^{-1}))^{-1} \leq t_{1}} \int_{a_0(\ell^*((4b)^{-1}))^{-1}}^{t_1}p_{d}(r,0)  |\varphi_{\alpha,\alpha+1} (s,r)| dr dy ds
\\
&=: II_{1}+II_{2}.
\end{align*}
Since $r\mapsto h(r)$ is decreasing, we see that $h((\ell^{-1}(a_0/r))^{-1})\ge h(4b)$ for $r\le a_0(\ell^*((4b)^{-1}))^{-1}$. Hence, by Proposition \ref{p:hku} (ii) and  Fubini's theorem
\begin{align*}
II_{1}&=\int_{\kappa(4b)}^{\infty} \int_{0}^{4b} \int_{a_0(\ell^*(\rho^{-1}))^{-1}}^{a_0(\ell^*((4b)^{-1}))^{-1}}(\ell^{-1}(a_{0}/r))^{d}e^{-crh((\ell^{-1}(a_0/r))^{-1})} s^{-\alpha-1}\rho^{d-1}dr d\rho ds\\
&\le\int_{\kappa(4b)}^{\infty} \int_{0}^{a_0(\ell^*((4b)^{-1}))^{-1}} \int_{0}^{(\ell^{-1}(a_0/r))^{-1}}(\ell^{-1}(a_{0}/r))^{d}e^{-crh(4b)} s^{-\alpha-1}\rho^{d-1}d\rho dr ds\\
&\le C\int_{\kappa(4b)}^{\infty} \int_{0}^{a_0(\ell^*((4b)^{-1}))^{-1}} e^{-crh(4b)} s^{-\alpha-1} dr ds\\
&\le C\int_{\kappa(4b)}^{\infty} \frac{1}{h(4b)} s^{-\alpha-1} ds\le C.
\end{align*}
Also for $II_{2}$, by using Proposition \ref{p:hku} (ii) and relation $K\leq h$ and $se^{-s} \leq c e^{-s/2}$ we have
\begin{align*}
II_{2}&\leq Cb^d\int_{\kappa(4b)}^{\infty}\int_{a_0(\ell^*((4b)^{-1}))^{-1}}^{t_{1}}p_{d}(a_{0}(\ell^{\ast}((4b)^{-1}))^{-1},0) s^{-\alpha-1}dr ds\\
&\le Cb^d\int_{\kappa(4b)}^{\infty}\int_{a_0(\ell^*((4b)^{-1}))^{-1}}^{t_{1}} b^{-d} (\ell^{\ast}((4b)^{-1}))^{-1} K(4b) e^{-c\frac{h(4b)}{\ell^{\ast}((4b)^{-1})}} s^{-\alpha-1}dr ds \\
&\le C \int_{\kappa(4b)}^{\infty} e^{-c\frac{h(4b)}{\ell^{\ast}((4b)^{-1})}} s^{-\alpha-1}dr ds\\
&\le C  h(4b) e^{-c\frac{h(4b)}{\ell^{\ast}((4b)^{-1})}} \leq C h(4b) e^{-c\frac{h(4b)}{\ell((4b)^{-1})}}\leq C.
\end{align*}
Thus, we obtain $II\le C$.

Now, we show \eqref{q1upper_case2_pf}. 
First, we see that
\begin{align*}
&\int_{\kappa(b)}^\infty \int_{b\leq |y| \leq h^{-1}(s^{-\alpha})} \int_{0}^{t_1} |D_{x}p_{d} (r,y)| |\varphi_{\alpha,\alpha+1} (s,r)| dr dy ds
\\
&\leq \int_{\kappa(b)}^\infty \int_{b\leq |y| \leq h^{-1}(s^{-\alpha})} III(s,y)  + \mathbf{1}_{a_0(\ell^*(|y|^{-1}))^{-1} \leq t_{1}} IV (s,y)  dy ds  ,
\end{align*}
where
\begin{align*}
III(s,y) &=  \int_{0}^{a_0(\ell^*(|y|^{-1}))^{-1}} |D_{x}p_{d} (r,y)| |\varphi_{\alpha,\alpha+1} (s,r)| dr,
\\
IV(s,y) &=  \int_{a_0(\ell^*(|y|^{-1}))^{-1}}^{t_{1}}|D_{x}p_{d} (r,y)| |\varphi_{\alpha,\alpha+1} (s,r)|dr .
\end{align*}
Like \eqref{eqn 05.30.18:42}, we have
\begin{align*}
&\int_{\kappa(b)}^\infty \int_{b\leq |y| \leq h^{-1}(s^{-\alpha})} III(s,y) dy ds
\\
&\leq C \int_{\kappa(b)}^\infty \int_{b}^{h^{-1}(s^{-\alpha})} \int_{0}^{\infty} \rho^{-2} e^{-crh(\rho)}s^{-\alpha-1} dr d\rho  ds 
\\
&\leq C \int_{b}^\infty \int_{(h(\rho))^{-1/\alpha}}^{\infty} \int_{0}^{\infty} \rho^{-2} e^{-crh(\rho)}s^{-\alpha-1} dr d\rho  ds 
\\
&\leq C \int_{b}^\infty \int_{(h(\rho))^{-1/\alpha}}^{\infty} \frac{1}{h(\rho)}  \rho^{-2} s^{-\alpha-1} ds d\rho  
\leq C \int_{b}^{\infty}  \rho^{-2} d\rho 
\leq Cb^{-1}.
\end{align*} 
Hence, we only need to control $IV$. Observe that by Remark \ref{r:hku_largetime} (i), Remark \ref{rmk 07.27.10:03} and \eqref{eqn 08.12.17:59}
\begin{align*}
&\int_{\kappa(b)}^\infty \int_{b\leq |y| \leq h^{-1}(s^{-\alpha})} IV(s,y) dy ds
\\
&\leq \int_{\kappa(b)}^\infty \int_{b\leq |y| \leq h^{-1}(s^{-\alpha})} \int_{a_0(\ell^*(|y|^{-1}))^{-1}}^{t_{1}} |y| p_{d+2}(r,0) s^{-\alpha-1}  dr dy ds
\\
&\leq C \int_{\kappa(b)}^\infty \int_{b\leq |y| \leq h^{-1}(s^{-\alpha})} \int_{a_0(\ell^*(|y|^{-1}))^{-1}}^{t_{1}} |y| p_{d+2}(a_{0}(\ell^{\ast}(|y|^{-1}))^{-1},0) s^{-\alpha-1}  dr dy ds
\\
&\leq C \int_{\kappa(b)}^\infty \int_{b\leq |y| \leq h^{-1}(s^{-\alpha})} \int_{a_0(\ell^*(|y|^{-1}))^{-1}}^{t_{1}} |y|^{-d-1}  e^{-c\frac{h(|y|)}{\ell^{\ast}(|y|^{-1})}}s^{-\alpha-1} dr dy ds
\\
& \leq C \int_{\kappa(b)}^\infty \int_{b}^{h^{-1}(s^{-\alpha})} \rho^{-2} \frac{h(\rho)}{h(\rho)} e^{-c\frac{h(\rho)}{\ell^{\ast}(\rho^{-1})}}s^{-\alpha-1} d\rho ds
\\
&\leq C \int_{\kappa(b)}^\infty \int_{b}^{h^{-1}(s^{-\alpha})} \rho^{-2} \frac{1}{h(\rho)}s^{-\alpha-1} d\rho ds
\\
& \leq C \int_{b}^{\infty}\int_{(h(\rho))^{-1/\alpha}}^{\infty} \rho^{-2} \frac{1}{h(\rho)}s^{-\alpha-1} ds d\rho 
\leq C b^{-1},
\end{align*}
Thus, we obtain \eqref{q1upper_case2_pf}. The lemma is proved.
\end{proof}

\end{document}